\documentclass[11pt,a4paper]{amsart}
\usepackage{a4wide}
\usepackage{amsmath,amssymb,url}

\date{\today}
\author[P. P. Pach]{P\'{e}ter P\'{a}l Pach}
\address{E\"otv\"os Lor\'and University,
          Department of Algebra and Number Theory,
         1117 Budapest,  P\'azm\'any P\'eter s\'et\'any 1/c, Hungary}
\email{ppp24@cs.elte.hu}

\author[M. Pinsker]{Michael Pinsker}
    \address{\'{E}quipe de Logique Math\'{e}matique\\ Universit\'{e} Diderot - Paris 7\\
	UFR de Math\'{e}matiques\\
	75205 Paris Cedex 13, France}
    \email{marula@gmx.at}
    \urladdr{http://dmg.tuwien.ac.at/pinsker/}
    \thanks{The second author is grateful for support through an APART-fellowship of the Austrian Academy of Sciences.}

\author[G. Pluh\'{a}r]{Gabriella Pluh\'{a}r}
\address{E\"otv\"os Lor\'and University,
          Department of Algebra and Number Theory,
         1117 Budapest,  P\'azm\'any P\'eter s\'et\'any 1/c, Hungary}
\email{plugab@cs.elte.hu}

\author[A. Pongr\'{a}cz]{Andr\'{a}s Pongr\'{a}cz}
\address{Central European University,
    Department of Mathematics and its Applications,
    1051 Budapest, N\'ador utca 9, Hungary}
\email{pongeee@cs.elte.hu}

\author[Cs. Szab\'{o}]{Csaba Szab\'{o}}
\address{E\"otv\"os Lor\'and University,
          Department of Algebra and Number Theory,
         1117 Budapest,  P\'azm\'any P\'eter s\'et\'any 1/c, Hungary}
\email{csaba@cs.elte.hu}

\title[Reducts of the random partial order]{Reducts of the random partial order}

\theoremstyle{plain}

    \newtheorem{thm}{Theorem}
    \newtheorem{theorem}[thm]{Theorem}
    \newtheorem{lem}[thm]{Lemma}
    
    \newtheorem{prop}[thm]{Proposition}
    \newtheorem{proposition}[thm]{Proposition}
    \newtheorem{cor}[thm]{Corollary}
    
    \newtheorem{fact}[thm]{Fact}

\theoremstyle{definition}

    \newtheorem{defn}[thm]{Definition}
    \newtheorem{definition}[thm]{Definition}
    \newtheorem{notation}[thm]{Notation}

\newcommand{\C}{\mathcal C}
\newcommand{\D}{\mathcal D}
\newcommand{\F}{\mathcal F}
\newcommand{\G}{\mathcal G}

\renewcommand{\P}{\mathbb P}
\renewcommand{\S}{\mathcal S}

\newcommand{\cl}[1]{\langle #1 \rangle}
\newcommand{\To}{\rightarrow}
\newcommand{\nin}{\notin}
\newcommand{\inv}{^{-1}}

\DeclareMathOperator{\Aut}{Aut}

\DeclareMathOperator{\Sym}{Sym}
\DeclareMathOperator{\id}{id}
\DeclareMathOperator{\tp}{tp}

\DeclareMathOperator{\Rev}{Rev}
\DeclareMathOperator{\R}{Rev}
\DeclareMathOperator{\T}{Turn}
\DeclareMathOperator{\M}{Max}
\DeclareMathOperator{\cyc}{Cycl}
\DeclareMathOperator{\pari}{Par}

\newcommand{\turn}{{\circlearrowright}}
\newcommand{\rev}{\updownarrow}
\newcommand{\Puv}{(\P,u,v)}
\newcommand{\oo}{\pari}
\newcommand{\oz}{\cyc}
\newcommand{\od}{\cyc'}
\newcommand{\tf}{\triangleleft_F}

\begin{document}
\begin{abstract}
	We determine, up to the equivalence of first-order interdefinability, all structures which are first-order definable in the random partial order. It turns out that these structures fall into precisely five equivalence classes. We achieve this result by showing that there exist exactly five closed permutation groups which contain the automorphism group of the random partial order, and thus expose all symmetries of this structure. Our classification lines up with previous similar classifications, such as the structures definable in the random graph or the order of the rationals; it also provides further evidence for a conjecture due to Simon Thomas which states that the number of structures definable in a homogeneous structure in a finite relational language is, up to first-order interdefinability, always finite. The method we employ is based on a Ramsey-theoretic analysis of functions acting on the random partial order, which allows us to find patterns in such functions and make them accessible to finite combinatorial arguments.
\end{abstract}
\maketitle
\section{Reducts of homogeneous structures}

The \emph{random partial order} $\P:=(P;\leq)$ is the unique countable partial order which is \emph{universal} in the sense that it contains all countable partial orders as induced suborders and which is \emph{homogeneous}, i.e., any isomorphism between two finite induced suborders of $\P$ extends to an automorphism of $\P$. Equivalently, $\P$ is the \emph{Fra\"{i}ss\'{e} limit} of the class of finite partial orders -- confer the textbook~\cite{Hodges}.

As the ``generic order'' representing all countable partial orders, the random partial order is of both theoretical and practical interest. The latter becomes in particular evident with the recent applications of homogeneous structures in theoretical computer science; see for example~\cite{BP-reductsRamsey,BodPin-Schaefer-both,tcsps-journal,MacphersonSurvey}. It is therefore tempting to classify all structures which are first-order definable in $\P$, i.e., all relational structures on domain $P$ all of whose relations can be defined from the relation $\leq$ by a first-order formula. Such structures have been called \emph{reducts} of $\P$ in the literature~\cite{RandomReducts,Thomas96}. It is the goal of the present paper to obtain such a classification \emph{up to first-order interdefinability}, that is, we consider two reducts $\Gamma, \Gamma'$ equivalent iff they are reducts of one another. We will show that up to this equivalence, there are precisely five reducts of $\P$. 

Our result lines up with a number of previous classifications of reducts of similar generic structures up to first-order interdefinability. The first non-trivial classification of this kind was obtained by Cameron~\cite{Cameron5} for the order of the rationals, i.e., the Fra\"{i}ss\'{e} limit of the class of finite linear orders; he showed that this order has five reducts up to first-order interdefinability. Thomas~\cite{RandomReducts} proved that the random graph has five reducts up to first-order interdefinability as well, and later generalized this result by showing that for all $k\geq 2$, the random hypergraph with $k$-hyperedges has $2^k+1$ reducts up to first-order interdefinability~\cite{Thomas96}. Junker and Ziegler~\cite{JunkerZiegler} showed that the structure $(\mathbb{Q};<,0)$, i.e., the order of the rationals with an additional constant symbol, has 116 reducts up to interdefinability. Further examples include the random $K_n$-free graph for all $n\geq 3$ (2 reducts, see~\cite{RandomReducts}), the random tournament (5 reducts, see~\cite{Bennett-thesis}), and the random $K_n$-free graph with a fixed constant (13 reducts if $n=3$ and 16 reducts if $n\geq 4$, see~\cite{Pon11}). A negative ``result'' is the random graph with a fixed constant, on which a subset of the authors of the present paper, together with another collaborator, gave up after having found 300 reducts. Obviously, the successful classifications have in common that the number of reducts is finite, and it is indeed an open conjecture of Thomas~\cite{RandomReducts} that all homogeneous structures in a finite relational language have only finitely many reducts up to first-order interdefinability.

The mentioned  classifications have all been obtained by means of the automorphism groups of the reducts, and we will proceed likewise in the present paper. It is clear that if $\Gamma$ is a reduct of a structure $\Delta$, then the automorphism group $\Aut(\Gamma)$ of $\Gamma$ is a permutation group containing $\Aut(\Delta)$, and also is a closed set with respect to the convergence topology on the space of all permutations on the domain of $\Delta$. If $\Delta$ is $\omega$-categorical, i.e., if $\Delta$ is up to isomorphism the only countable model of its first-order theory, then it follows from the theorem of Ryll-Nardzewski, Engeler and Svenonius (confer~\cite{Hodges}) that the converse is true as well: the closed permutation groups acting on the domain of $\Delta$ and containing $\Aut(\Delta)$ are precisely the automorphism groups of reducts of $\Delta$; moreover, two reducts have equal automorphism groups if and only if they are first-order interdefinable. Since homogeneous structures in a finite language are $\omega$-categorical, it is enough for us to determine all closed permutation groups that contain $\Aut(\P)$ in order to obtain our classification.

The fact that the reducts of an $\omega$-categorical structure $\Delta$ correspond to the closed permutation groups containing $\Aut(\Delta)$ not only yields a method for classifying these reducts, but also a meaningful interpretation of such classifications: for just like $\Aut(\Delta)$ is the group of all symmetries of $\Delta$, the closed permutation groups containing $\Aut(\Delta)$ stand for all symmetries of $\Delta$ if we are willing to give up some of the structure of $\Delta$. As for an example, it is obvious that turning the random partial order upside down, one obtains again a random partial order; this symmetry is reflected by one of the closed groups containing $\Aut(\P)$, namely the group of all automorphisms and antiautomorphisms of $\P$. It will follow from our classification that $\P$ has only one more symmetry of this kind -- this second symmetry is much less obvious, and so we argue that the classification of the reducts of $\P$, or indeed of any $\omega$-categorical structure, is much more than a mere sportive challenge -- it is an essential part of understanding the structure itself.

Our approach to investigating the closed groups containing $\Aut(\P)$ is based on a Ramsey-theoretic analysis of functions, and in particular permutations, on the domain $P$ of $\P=(P;\leq)$; this allows us to find patterns of regular behaviour with respect to the structure $\P$ in  any arbitrary function acting on $P$. The method as we use it has been developed in~\cite{BodPinTsa, BodPin-Schaefer-both,RandomMinOps,BP-reductsRamsey} and is a general powerful technique for dealing with functions on \emph{ordered homogeneous Ramsey structures} in a finite language. But while this machinery has previously been used, for example, to re-derive and extend Thomas' classification of the reducts of the random graph, it is only in the present paper (and, at the same time, in~\cite{Pon11} for the reducts of $K_n$-free graphs with a constant)  that it is applied to obtain a new full classification of reducts of a homogeneous structure up to first-order interdefinability.

Before stating our result, we remark that finer classifications of reducts of homogeneous structures, for example up to existential, existential positive, or primitive positive interdefinability,  have also been considered in the literature, in particular in applications -- see \cite{BodChenPinsker, BodPinTsa, RandomMinOps,BP-reductsRamsey}.

\section{The reducts of the random partial order}

\subsection{The group formulation}

In a first formulation of our result, we will list the closed groups containing $\Aut(\P)$ by means of sets of permutations generating them: we say that a set $\S$ of permutations  on $P$ \emph{generates} a permutation $\alpha$ on $P$ iff $\alpha$ is an element of the smallest closed permutation group $\cl{\S}$ that contains $\S$. Equivalently, writing $\id$ for the identity function on $P$, for every finite set $A\subseteq P$ there exist $n\geq 0$, $\beta_1,\ldots,\beta_n\in \S$, and $i_1,\ldots,i_n\in\{1,-1\}$ such that $\beta_1^{i_1}\circ\cdots\circ \beta_n^{i_n}\circ\id$ agrees with $\alpha$ on $A$. We also say that a permutation $\beta$ generates $\alpha$ iff $\{\beta\}$ generates $\alpha$.

If for $x,y\in P$ we define $x \geq y$ iff $y\leq x$, then the structure $(P;\geq)$ is isomorphic to $\P$ -- it is, for example, easy to verify that it contains all finite partial orders and that it is homogeneous. Hence, there exists an isomorphism between the two structures, and we fix one such isomorphism $\rev: P\To P$; so the function $\updownarrow$ simply reverses the order $\leq$ on $P$. It is easy to see that any two isomorphisms of this kind generate one another, and the exact choice of the permutation is thus irrelevant for our purposes.

The class $\C$ of all finite structures of the form $(A;\leq',F')$, where $\leq'$ is a partial order on $A$, and $F'\subseteq A$ is an upward closed set with respect to $\leq'$, is an \emph{amalgamation class} in the sense of~\cite{Hodges}. Hence, it has a Fra\"{i}ss\'{e} limit; that is, there exists an up to isomorphism unique countable structure which is homogeneous and whose \emph{age}, i.e., the set of finite structures isomorphic with one of its induced substructures, equals $\C$. The partial order of this limit is just the random partial order, and thus we can write $(P;\leq,F)$ for this structure, where $F\subseteq P$ is an upward closed set with respect to $\leq$. By homogeneity and universality of $(P;\leq,F)$, $F$ is even a \emph{filter}, i.e., any two elements of $F$ have a lower bound in $F$. We call $(P;\leq,F)$ the \emph{random partial order with a random filter}, and any filter $W\subseteq P$ with the property that $(P;\leq, W)$ is isomorphic with $(P;\leq,F)$ \emph{random}.

Let $F\subseteq P$ be a random filter, and let $I:=P\setminus F$. Then $I$ is downward closed, and in fact an \emph{ideal}, i.e., any two elements of $I$ have an upper bound in $I$. Define a partial order $\trianglelefteq_F$ on $P$ by setting 
\begin{align*}
	x \trianglelefteq_F y\; \leftrightarrow\; & x,y\in F\text{ and } x\leq y, \text{ or }\\
	& x,y\in I \text{ and } x\leq y, \text{ or }\\ 
	&x\in F\wedge y\in I\text{ and } y\nleq x, 
\end{align*}
where $a\nleq b$ is short for $\neg(a\leq b)$. It is easy to see that $(P;\trianglelefteq_F)$ is indeed a partial order, and we will verify in the next section that $(P;\trianglelefteq_F)$ and $\P$ are isomorphic. Pick an isomorphism $\turn_F\colon (P;\trianglelefteq_F)\To\P$. Then for $x,y\in F$, we have $f(x)\leq f(y)$ if and only if $x\leq y$, and likewise for $x,y\in I$; if $x\in F$ and $y\in I$, then $f(x)\leq f(y)$ if and only if $y\nleq x$; and moreover, $f(x)\ngeq f(y)$ for all $x\in F$ and $y\in I$. 
It is not hard to see that any two permutations obtained this way generate one another, even if they were defined by different random filters. We therefore also write $\turn$ for any $\turn_F$ when the filter $F$ is not of particular interest.

\begin{theorem}\label{thm:groups}
	The following five groups are precisely the closed permutation groups on $P$ which contain $\Aut(\P)$.
	\begin{enumerate}
	\item $\Aut(\P)$;
	\item $\R:=\cl{\{\rev\}\cup\Aut(\P)}$;
	\item $\T:=\cl{\{\turn\}\cup\Aut(\P)}$;
	\item $\M:=\cl{\{\rev,\turn\}\cup\Aut(\P)}$;
	\item The full symmetric group $\Sym_P$ of all permutations on $P$.
	\end{enumerate}
\end{theorem}

As a consequence, the only symmetries of $\P$ in the sense mentioned above are turning it upside down, and ``turning'' it around a random filter $F$ via the function $\turn_F$. These symmetries suggest the investigation of the corresponding operations on finite posets (essentially, the restrictions of $\rev$ and $\turn_F$ to finite substructures of $\P$). While $\rev$ for finite posets is, of course, combinatorially not very exciting, the study of ``turns'' of finite posets seems to be quite worthwhile -- we refer to the companion paper~\cite{newTransformation}.

We will also obtain explicit descriptions of the elements of the groups in Theorem~\ref{thm:groups}. 
Clearly, the group $\R$ contains exactly the automorphisms of $\P$ and the isomorphisms between $\P$ and $(P;\geq)$. We will show that $\T$ consists precisely of what we will call \emph{rotations} in Definition~\ref{defn:rotations} -- these are functions of slightly more general form that the functions $\turn_F$. Moreover, $\M$ turns out to be simply the union of $\R$, $\T$, and the set of all functions of the form $\rev\circ f$, where $f$ is a rotation. 

\subsection{The reduct formulation}

We now turn to the relational formulation of our result; that is, we will specify five reducts of $\P$ such that any reduct of $\P$ is first-order interdefinable with one of the reducts of our list.

 Define a binary relation $\bot$ on $P$ by $\bot:=\{(x,y)\in P^2\;|\; x\nleq y\wedge y\nleq x \}$. We call the relation the \emph{incomparability relation}, and refer to elements $x,y\in P$ as \emph{incomparable} iff $(x,y)$ is an element of $\bot$; in that case, we also write $x\bot y$. Elements $x,y\in P$ are \emph{comparable} iff they are not incomparable.
 
 For $x,y\in P$, write $x<y$ iff $x\leq y$ and $x\neq y$. Now define a ternary relation $\cyc$ on $P$ by
 \begin{align*}
 	\cyc:=\{(x,y,z)\in P^3\;|\;& (x<y<z)\vee (y<z<x)\vee (z<x<y)\vee\\ &(x<y\wedge x\bot z\wedge y\bot z) \vee\\& (y<z\wedge y\bot x\wedge z\bot x) \vee\\& (z<x\wedge z\bot y\wedge x\bot y)\}.
 \end{align*}

Finally, define a ternary relation $\pari$ on $P$ by
\begin{align*}
	\pari:=\{(x,y,z)\in P^3\;|&\; x,y,z\text{ are distinct  and the number of }\\
	&\text{ 2-element subsets  of incomparable elements of } \{x,y,z\} \text{ is odd}\}.
\end{align*}

\begin{theorem}\label{thm:reducts}
	Let $\Gamma$ be a reduct of $\P$. Then $\Gamma$ is first-order interdefinable with precisely one of the following structures.
	\begin{enumerate}
	\item $\P=(P;\leq)$;
	\item $(P;\bot)$;
	\item $(P;\cyc)$;
	\item $(P;\pari)$;
	\item $(P;=)$.
	\end{enumerate}
	Moreover, for $1\leq x\leq 5$, $\Gamma$ is first-order interdefinable with structure (x) if and only if $\Aut(\Gamma)$ equals group number (x) in Theorem~\ref{thm:groups}.
\end{theorem}

\section{Random filters and the extension property}

Before turning to the main proof of our theorems, we verify the existence of the permutation $\turn_F$. That is, we must show that if $F\subseteq P$ is a random filter, then $(P;\triangleleft_F)$ and $\P$ are isomorphic. The easiest way to see this is by checking that $(P;\trianglelefteq_F)$ satisfies the following \emph{extension property}, which determines $\P$ up to isomorphism and which we will use throughout the paper: 
for any finite set $S=\{s_1,\ldots,s_k\}\subseteq P$ and any partial order with domain $\{y\}\cup S$ extending the order induced by $\P$ on $S$, there exists $x\in P$ such that the assignment from $\{x\}\cup S$ to $\{y\}\cup S$ which sends $x$ to $y$ and leaves all elements of $S$ fixed is an isomorphism. In logic terminology, the extension property says that if we fix any finite set of elements $s_1,\ldots,s_k\in P$, and express properties of another imaginary element $x$ by means of a quantifier-free $\{\leq\}$-formula with one free variable using parameters $s_1,\ldots,s_k$, then an element enjoying these properties actually exists in $\P$ unless the properties are inconsistent with the theory of partial orders.

\begin{prop}
Let $F\subseteq P$ be a random filter of $\P$. Then $(P;\triangleleft_F)$ satisfies the extension property. Consequently, $(P;\triangleleft_F)$ and $\P$ are isomorphic and $\turn_F$ exists. 
\end{prop}
\begin{proof}
Let $s_1,\ldots,s_k\in P$ and an extension of the order induced by $\tf$ on $S=\{s_1,\ldots,s_k\}$ by an element $y$ outside $S$ be given. We will denote the order on $T:=S\cup\{y\}$ by $\triangleleft_F$ as well. Let $I:=P\setminus F$ be the ideal in $\P$ corresponding to the filter $F$, and write $S$ as a disjoint union $S_F\cup S_I$, where $S_F:=S\cap F$, and $S_{I}:=S\cap I$. Now suppose that there exist $a\in S_{I}$ and $b\in S_{F}$ such that $a\triangleleft_F y \triangleleft_F b$. Then $a\triangleleft_F b$, which is impossible by the definition of $\triangleleft_F$, since $a\in {I}$ and $b\in F$. Hence, assume without loss of generality that we do not have $y\triangleleft_F b$ for any $b\in S_F$. Then $W:=S_I\cup\{y\}$ is  upward closed and $S_F$ downward closed in $(T;\tf)$. Now define an order $\leq_W$ on $T$ by setting 
$$
	u \leq_W v\; \leftrightarrow (u,v\in W\wedge u\tf v)\vee (u,v\in T\setminus W\wedge u\tf v)\vee (u\in W\wedge v\in T\setminus W\wedge \neg(v\tf u)).
$$ 
Note that this defines $\leq_W$ from $\tf$ in precisely the same way as $\tf$ is defined (though on $\P$) from $\leq$. Hence, $\leq_W$ is a partial order on $T$, and for $u,v\in S$ we have $u\leq_V v$ if and only if $u\leq v$ in $\P$. Now the downward closed set $S_F$ in $(T;\tf)$ is an upward closed set in $(T;\leq_W)$. Hence, the structure $(T;\leq_W,S_F)$ has an embedding $\xi$ into the universal object $(P;\leq,F)$. Since $\leq_W$ agrees with $\leq$ on $S$, and by homogeneity, we may assume that $\xi$ is the identity on $S$. Set $x:=\xi(y)$. We leave the straightforward verification of the fact that 
the assignment from $\{x\}\cup S$ to $\{y\}\cup S$ which sends $x$ to $y$ and leaves all elements of $S$ fixed is an isomorphism from $(\{x\}\cup S;\tf)$ onto $(T;\tf)$ to the reader.
\end{proof}

Let us remark that the ideal $I=P\setminus F$ corresponding to a random filter $F$ on $\P$ is random in the analogous sense for ideals. Moreover, under $\turn_F$ the random filter $F$ is sent to a random ideal, and vice-versa. One could thus assume that the image of $F$ under $\turn_F$ equals $I$, in which case $\turn_F$ becomes, similarly to $\rev$, its own ``almost'' inverse in the sense that applying it twice yields an automorphism of $\P$. By adjusting it with such an automorphism, one could even assume that $\turn_F=\turn_F\inv$.

\section{Ramsey theory: canonizing functions}

Our combinatorial method for proving Theorem~\ref{thm:groups} is to apply Ramsey theory in order to find patterns of regular behaviour in arbitrary functions on $\P$, and follows~\cite{BodPinTsa, BodPin-Schaefer-both,RandomMinOps,BP-reductsRamsey}. We make this more precise.

\begin{definition}\label{defn:type}
    Let $\Delta$ be a structure. The \emph{type} $\tp(a)$ of an $n$-tuple $a$ of elements in $\Delta$ is the set of first-order formulas with
     free variables $x_1,\dots,x_n$ that hold for $a$ in $\Delta$.
\end{definition}

%We bring to the reader's attention the well-known fact that in homogeneous structures in a finite language, in particular in %the random partial order, two $n$-tuples $a,b$ have the same type if and only if $a$ can be sent to $b$ by an %automorphism.

\begin{definition}\label{defn:behaviour}
    Let $\Delta, \Lambda$ be structures. A \emph{type condition} between $\Delta$ and $\Lambda$ is a pair $(t,s)$, where $t$ is a type of an $n$-tuple in $\Delta$, and $s$ is a type of an $n$-tuple in $\Lambda$, for some $n\geq 1$. 
    
      A function $f:\Delta\To \Lambda$ \emph{satisfies} a type condition $(t,s)$ between $\Delta$ and $\Lambda$ iff for all $n$-tuples $a=(a_1,\ldots,a_n)$ of elements of $\Delta$ with $\tp(a)=t$ the $n$-tuple $f(a):=(f(a_1),\ldots,f(a_n))$ has type $s$ in $\Lambda$. A \emph{behaviour} is a set of type conditions between structures $\Delta$ and $\Lambda$. A function from $\Delta$ to $\Lambda$ \emph{has behaviour $B$} iff it satisfies all the type conditions of $B$.
\end{definition}

\begin{definition}\label{defn:canonical}
    Let $\Delta, \Lambda$ be structures. A function $f: \Delta \To \Lambda$ is \emph{canonical} iff for all types $t$ of $n$-tuples in $\Delta$ there exists a type $s$ of an $n$-tuple in $\Lambda$ such that $f$ satisfies the type condition $(t,s)$. In other words, $n$-tuples of equal type in $\Delta$ are sent to $n$-tuples of equal type in $\Lambda$ under $f$, for all $n\geq 1$.
\end{definition}

We remark that since $\P$ is homogeneous, every first-order formula is over $\P$ equivalent to a quantifier-free formula, and so the type of an $n$-tuple $a$ in $\P$ is determined by which of its elements are equal and between which elements the relation $\leq$ holds. In particular, the type of $a$ only depends on its binary subtypes, i.e., the types of the pairs $(a_i, a_j)$, where $1\leq i,j\leq n$. Therefore, a function $f: \P\To \P$ is canonical iff it satisfies the condition of the definition for types of 2-tuples.

Roughly, our strategy is to make the functions we work with canonical, and thus easier to handle. To achieve this, we first enrich the structure $\P$ by a linear order in order to improve its combinatorial properties, as follows. We do not give the -- in some cases fairly technical -- definitions of all notions in this discourse, as they will not be needed later on; in any case, Proposition~\ref{prop:orderedRandomCanonical} that follows is used as a black box for this paper, and the reader interested in its proof is referred to~\cite{BodPinTsa}.
The class $\D$ of all finite structures $(A;\leq',\prec')$ with two binary relations $\leq'$ and $\prec'$, where $\leq'$ is a partial order and $\prec'$ is a total order extending $\leq'$, is an amalgamation class, and moreover a \emph{Ramsey class} (see for example~\cite[Theorem~1~(1)]{Sokic}). By the first property, it has a Fra\"{i}ss\'{e} limit.  Checking the extension property, one sees that the partial order of this limit is just the random partial order, and by uniqueness of the dense linear order without endpoints its total order is isomorphic to the order of the rationals. Hence, there exists a linear order $\prec$ on $P$ which is isomorphic to the order of the rationals, which extends $\leq$, and such that the structure $\P^+:=(P;\leq,\prec)$ is precisely the Fra\"{i}ss\'{e} limit of the class $\D$. So $\P^+$ is a homogeneous structure in a finite language which has a  linear order among its relations and which is \emph{Ramsey}, i.e. its age, which equals the class $\D$, is a Ramsey class. The following proposition is then  a consequence of the results in~\cite{BodPinTsa,BP-reductsRamsey} about such structures. To state it, let us extend the notion ``generates'' to non-permutations: for a set of functions $\F\subseteq P^P$ and $f\in P^P$, we say that $f$ is 
\emph{M-generated} by $\F$ iff it is contained in the smallest transformation monoid on $P$ which contains $\F$ and which is a closed set in the convergence topology on $P^P$. In other words, $f$ is M-generated by $\F$ iff for all finite $A\subseteq P$ there exist $n\geq 0$ and $f_1,\ldots,f_n\in\F$ such that $f_1\circ\cdots\circ f_n\circ \id$ agrees with $f$ on $A$. For a structure $\Delta$ and elements $c_1,\ldots,c_n$ of $\Delta$, we write $(\Delta,c_1,\ldots,c_n)$ for the structure obtained by adding the constant symbols $c_1,\ldots,c_n$ to $\Delta$.

\begin{proposition}\label{prop:orderedRandomCanonical}
	Let $f: P\To P$ be a function, and let $c_1,\ldots,c_n,d_1,\ldots,d_m\in P$. Then $\{f\}\cup\Aut(\P^+)$ M-generates a function which is canonical as a function from $(\P^+,c_1,\ldots,c_n)$ to $(\P^+,d_1,\ldots,d_m)$, and which is identical with $f$ on $\{c_1,\ldots,c_n\}$.
\end{proposition}

Any canonical function $g$ from $(\P^+,c_1,\ldots,c_n)$ to $(\P^+,d_1,\ldots,d_m)$ defines a function from the set $T$ of types of pairs of distinct elements in $(\P^+,c_1,\ldots,c_n)$ to the set $S$ of such types in $(\P^+,d_1,\ldots,d_m)$ -- this ``type function'' simply assigns to every element $t$ of  $T$ the type $s$ in $S$ for which the type condition $(t,s)$ is satisfied by $g$. Already when $n=m=0$, i.e., there are no constants added to $\P^+$, then $|T|=|S|=4$, so in theory there are $4^4$ such type functions. The following lemma states which of them actually occur. 

\begin{lem}\label{lem:behaviours}
	Let $g:\P^+\to\P^+$ be canonical and injective. Then it has one of the following behaviours.
	\begin{itemize}
		\item[(i a)] $g$ behaves like $\id$, i.e., it preserves $\leq$ and $\bot$ (and hence also $\prec$);		\item[(i b)] $g$ behaves like $\rev$, i.e., it reverses $\leq$ and preserves $\bot$ (and hence reverses $\prec$);
		\item[(ii a)] $g$ sends $P$ order preservingly onto a chain with respect to $\leq$  (and hence preserves $\prec$);
		\item[(ii b)] $g$ sends $P$ order reversingly onto a chain with respect to $\leq$  (and hence reverses $\prec$);

		\item[(iii a)] $g$ sends $P$ onto an antichain with respect to $\leq$ and preserves $\prec$;
		\item[(iii b)] $g$ sends $P$ onto an antichain with respect to $\leq$ and reverses $\prec$.
	\end{itemize}
\end{lem}
\begin{proof}
We first prove that $g$ either preserves or reverses the order $\prec$.

Suppose there exist $a,b\in P$ with $a\prec b$ such that $g(a)\prec g(b)$. Assume first that $a\leq b$. Then $g(c)\prec g(d)$ for all $c,d\in P$ with $c\prec d$ and $c\leq d$ because $g$ is canonical. Now using the universality of $\P^+$, pick $u,v,w\in P$ with $u\prec v\prec w$, $u\leq w$, $u\bot v$, and $v\bot w$. Then $g(u)\leq g(w)$ by our observation above. If $g(v)\prec g(u)$, then also $g(w)\prec g(v)$ as $g$ is canonical, and hence $g(w)\prec g(u)$, a contradiction. Hence, $g(u)\prec g(v)$, and so $g(c)\prec g(d)$ for all $c, d\in P$ with $c\prec d$, so $g$ preserves $\prec$. Now suppose that $a\bot b$. Then $g(c)\prec g(d)$ for all $c,d\in P$ with $c\prec d$ and $c\bot d$, because $g$ is canonical. Pick  $u,v,w\in P$ as before. This time, $g(u)\prec g(v)\prec g(w)$, and hence $g(u)\prec g(w)$. Therefore, $g(c)\prec g(d)$ for all $c, d\in P$ with $c\prec d$, so $g$ again preserves $\prec$. 

By the dual argument, the existence of $a,b\in P$ with $a\prec b$ such that $g(b)\prec g(a)$ implies that $g$ reverses $\prec$.

We next show that if $g$ preserves $\prec$, then one of the situations (i a), (ii a), (iii a) occurs; then by duality, if $g$ reverses $\prec$, one of (i b), (ii b), (iii b) hold. We distinguish two cases. 

Suppose first that $g(a)\bot g(b)$ for all $a,b\in P$ with $a\leq b$. Let  $c,d,e\in P$ such that $c\prec d\prec e$, $c\bot d$, $c\leq e$, and $e\bot d$. If $g(c)$ and $g(d)$ were comparable, then $g(c)\leq g(d)$ since $g(c)\prec g(d)$, and likewise $g(d)\leq g(e)$, so that $g(c)\leq g(d)$, a contradiction. Hence, $g(c)$ and $g(d)$ are incomparable, and so, since $g$ is canonical, (iii a) holds.

Assume now that $g(a)\leq g(b)$ for all  $a,b\in P$ with $a\leq b$. If $g(c)\leq g(d)$ also for all $c,d\in P$ with $c\bot d$ and $c\prec d$, then clearly we have situation (ii a). Otherwise, $g(c)\bot g(d)$ for all $c,d\in P$ with $c\bot d$ and $c\prec d$, and we have case (i a).

Since one of these two situations must be the case, we are done.
\end{proof}

When applying Proposition~\ref{prop:orderedRandomCanonical}, we will be able to ignore most of the possible behaviours of canonical functions as a consequence of the following lemma.

\begin{lem}\label{lem:chainantichain}
	Let $\G\supseteq \Aut(\P)$ be a closed group such that for all finite $A \subseteq P$ there is a function M-generated by $\G$ which sends $A$ to a chain or an antichain. Then $\G=\Sym_P$.
\end{lem}

\begin{proof}
	Suppose first that for all finite $A \subseteq P$ there is a function M-generated by $\G$ which sends $A$ to an antichain. Let $s,t$ be injective $n$-tuples of elements in $P$, for some $n\geq 1$. Let $g:P\To P$ and  $h: P\To P$ be functions M-generated by $\G$ such that $g(s)$ (the $n$-tuple obtained by applying $g$ to every component of $s$) and $h(t)$ induce antichains in $\P$. By the homogeneity of $\P$, there exists an automorphism $\alpha\in\Aut(\P)$ such that $\alpha(g(s))=h(t)$. Also, since $\G$ contains the inverse of all of its functions, there exists a function $p:P\To P$ M-generated by $\G$ such that $p(h(t))=t$, and hence $p(\alpha(g(s)))=t$. Since $p\circ\alpha\circ g$ is M-generated by $\G$, there exists $\beta\in\G$ which agrees with this function on $s$. Hence, $\beta(s)=t$, proving that $\G$ is $n$-transitive for all $n\geq 1$, and so $\G=\Sym_P$.
	
Now suppose that for all finite $A \subseteq P$ there is a function M-generated by $\G$ which sends $A$ to a chain. Let any finite $A\subseteq P$ be given, and let $B\subseteq P$ be so that $|B|=|A|$ and such that $B$ induces an independent set in $\P$. Let $g:P\To P$ and  $h: P\To P$ be functions M-generated by $\G$ such that $g[A]$ and $h[B]$ induce chains in $\P$. There exists $\alpha\in\Aut(\P)$ such that $\alpha[g[A]]=h[B]$. Let $p:P\To P$ be a function generated by $\G$ such that $p[h[B]]=B$. Then $p[\alpha[g[A]]]=B$, and hence we are back in the preceding case.

Finally, observe that one of the two cases must occur: for otherwise, there exist finite $A_1, A_2\subseteq P$ such that $A_1$ cannot be set to an antichain, and $A_2$ cannot be sent to a chain by any function which is M-generated by $\G$. But then $A_1\cup A_2$ can neither be sent to a chain nor to an antichain by any such function, a contradiction.
\end{proof}

\begin{lem}\label{lemma:ruleOutChainAntichain}
	Let $\G\supseteq \Aut(\P)$ be a closed group which M-generates a canonical function of behaviour (ii a), (ii b), (iii a) or (iii b) in Lemma~\ref{lem:behaviours}. Then $\G=\Sym_P$.
\end{lem}
\begin{proof}
	This is a direct consequence of Lemma~\ref{lem:chainantichain}.
\end{proof}

Having enriched $\P$ with the linear order $\prec$ and taken advantage of Proposition~\ref{prop:orderedRandomCanonical}, we pass to a suitable substructure of $(\P^+,c_1,\ldots,c_n)$ in order to get rid of $\prec$ -- this substructure will be called a \emph{$\prec$-clean skeleton}. Before giving the exact definition, we need more notions and notation concerning the definable subsets of $(\P,c_1,\ldots,c_n)$ and of $(\P^+,c_1,\ldots,c_n)$.

\begin{definition}
	Let $\G$ be a permutation group acting on a set $D$. Then for $n\geq 1$ and $a=(a_1,\ldots,a_n)\in D^n$, the set
	$$
		\{(\alpha(a_1),\ldots,\alpha(a_n)):\alpha\in\G\}\subseteq D^n
	$$
	is called an $n$-orbit of $\G$. The $1$-orbits are just called \emph{orbits}. If $\Delta$ is a structure, then the $n$-orbits of $\Delta$ are defined as the $n$-orbits of $\Aut(\Delta)$. 
\end{definition}

By the theorem of Ryll-Nardzewski, Engeler and Svenonius, two $n$-tuples in an $\omega$-categorical structure  belong to the same $n$-orbit if and only if they have the same type; in particular, this is true in the structures $(\P,c_1,\ldots,c_n)$ and $(\P^+,c_1,\ldots,c_n)$.

\begin{notation}\label{nota:orbits}
Let $c_1,\ldots,c_n\in P$. For $R_1,\ldots,R_n\in \{=,<, \perp, >\}$ and $S_1,\ldots,S_n\in \{ \prec,\succ\}$, we set
$$
	X_{R_1,\ldots, R_n}:=\{x\in P: c_1R_1 x\ \wedge\ \cdots\ \wedge\ c_nR_n x\}
$$
and
$$
	X_{R_1,\ldots, R_n}^{S_1,\ldots,S_n}:=\{x\in P: (c_1R_1 x\ \wedge c_1 S_1 x)\ \wedge\ \cdots\ \wedge\ (c_nR_n x\wedge x_n S_n x)\}.
$$
The constants $c_1,\ldots,c_n$ are not specified in the notation, but will always be clear from the context.
\end{notation}

The following is well-known and easy to verify using the homogeneity and universality of $\P$ and $\P^+$, and in particular the fact that first-order formulas over these structures are equivalent to quantifier-free formulas.

\begin{fact}\label{fact:orbits}
Let $c_1,\ldots,c_n\in P$.	
\begin{itemize}
	\item The sets $X_{R_1,\ldots, R_n}$ are either empty, or equal to $\{c_i\}$ for some $1\leq i\leq n$, or infinite and induce $\P$. The orbits of $(\P,c_1,\ldots,c_n)$ are precisely the non-empty sets  of this form.
\item The sets 
$X_{R_1,\ldots, R_n}^{S_1,\ldots,S_n}$ are either empty, or equal to $\{c_i\}$ for some $1\leq i\leq n$, or infinite and induce $\P^+$. The orbits of $(\P^+,c_1,\ldots,c_n)$ are precisely the non-empty sets of this form.
\end{itemize}
\end{fact}

\begin{defn}\label{defn:skeleton}
	Let $\Delta$ be a structure on domain $D$. A subset $S$ of $D$ is called a \emph{skeleton} of $\Delta$ iff it induces a substructure of $\Delta$ which is isomorphic to $\Delta$. Now let $\sqsubset$ be a linear order on $D$. Then a skeleton $S$ is called \emph{$\sqsubset$-clean} iff whenever $a=(a_1,a_2),b=(b_1,b_2)\in S^2$ have the same type in $\Delta$, then either $a,b$ or $a, \tilde b:=(b_2,b_1)$ have the same type in $(\Delta,\sqsubset)$.
\end{defn}

In this paper, we only need a $\prec$-clean skeleton of $(\P,c_1,\ldots,c_n)$, but we stated Definition~\ref{defn:skeleton} generally since we believe it could be useful in other situations where a homogeneous structure is extended by a linear order with the goal of making it Ramsey.

\begin{lem}\label{lem:existenceSkeleton}
	Let $c_1,\ldots,c_n\in P$. Then $(\P,c_1,\ldots,c_n)$ has a skeleton which is $\prec$-clean.
\end{lem}
\begin{proof}
	Let $O_1,\ldots,O_k$ be the orbits of $(\P,c_1,\ldots,c_n)$, and pick one representative element $r_i$ of each orbit $O_i$. By relabelling the orbits, we may assume that $r_1\prec\cdots\prec r_k$; pick an additional $r_0\in P$ with $r_0\prec r_1$. 
	Now for all $1\leq j\leq k$ for which $O_j$ is infinite set
	$$
	S_j:= \{s\in O_j\; |\; r_{j-1}\prec s\prec r_j\}. 
	$$
	Let $S$ be the union of all the $S_j$ with $\{c_1,\ldots,c_n\}$. To see that $S$ is a skeleton, it suffices to verify the extension property for $(S;\leq)$. Let $U=\{u_1,\ldots,u_l\}\subseteq S$ induce a finite substructure of $(S;\leq)$, and let $U\cup\{y\}$ be an extension of $U$ by an element $y\nin U$. We may assume that $U$ contains $\{c_1,\ldots,c_n\}$. By the extension property for $\P$, we may assume that $y$ is an element of this structure, and so $y\in O_j$ for some $1\leq j\leq k$. Since $y\nin \{c_1,\ldots,c_n\}$, $O_j$ is infinite. We claim there exists $x\in S_j$ such that $x,y$ have the same type in $(\P,u_1,\ldots,u_l)$ -- then picking any such $x$ yields the desired extension. Otherwise, let $\phi(z)$ be a conjunction as in the first part of Notation~\ref{nota:orbits} which defines $O_j$ in $(\P,u_1,\ldots,u_l)$, i.e., $\phi(z)$ is the conjunction of all atomic formulas with one free variable $z$ satisfied by $y$ in this structure. By our assumption, $\phi(z)$ implies $z\nin S_j$, so it implies $r_{j-1}\not \prec z\vee  z\not\prec r_j$. By the universality and homogeneity of $\P^+$, this is only possible if $\phi(z)$ implies $z\leq r_{j-1}\vee r_j\leq z$ in $\P$, which can only happen if there exists $1\leq i\leq l$ such that $y<u_i\leq r_{j-1}$ or $r_j\leq u_i<y$. Consider the second case; the first case is isomorphic. By the definition of $S$, we conclude $u_i\in S_p$ for some $p\neq j$, and so the orbits of $u_i$ and $r_j$ in $(\P,c_1,\ldots,c_n)$ are distinct. Therefore,  there exists $1\leq m\leq n$ such that either $u_i\geq  c_m$ and $r_j\not \geq c_m$, or $r_j< c_m$ and $u_i \not <  c_m$. In  the first case we infer $y> c_m$, contradicting the fact that $y$ and $r_j$ have the same type in $(\P,c_1,\ldots,c_n)$. In the second case it follows that $y\not < c_m$, yielding the same contradiction.
 
We show that $S$ is $\prec$-clean. Let $a=(a_1,a_2),b=(b_1,b_2)\in S^2$ have the same type in $(\P,c_1,\ldots,c_n)$. Then there exist $1\leq i,j\leq k$ such that $a_1,b_1\in O_i$ and $a_2,b_2\in O_j$. Suppose $i=j$. If $a_1,a_2$ are comparable, say $a_1\leq a_2$, then $b_1\leq b_2$, $a_1\prec a_2$, and $b_1\prec b_2$, and we are done. If $a_1\bot a_2$, then  $b_1\bot b_2$ and so either $a,b$ or $a,\tilde b=(b_2,b_1)$ have the same type in $(\P^+,c_1,\ldots,c_n)$. Now suppose $i\neq j$, say $i<j$. Then $a_1\prec a_2$ and  $b_1\prec b_2$, and so $a,b$ have the same type in $(\P^+,c_1,\ldots,c_n)$.
 
\end{proof}

\begin{lem}\label{lem:interface}
	Let $f: P\To P$ be a permutation, and let $c_1,\ldots,c_n,d_1,\ldots,d_m\in P$. Then $\{f,f\inv \}\cup\Aut(\P)$ M-generates a function $g:P\To P$ with the following properties.
	\begin{itemize}
		\item $g$ agrees with $f$ on $\{c_1,\ldots,c_n\}$;
		\item $g$ is canonical as a function from $(\P,c_1,\ldots,c_n)$ to $(\P,d_1,\ldots,d_m)$.
%		\item on each infinite orbit of $(\P,c_1,\ldots,c_n)$, $g$ behaves like $\id$ or like $\rev$.
	\end{itemize}
\end{lem}

\begin{proof}
	  Let $h$ be the function guaranteed by Proposition~\ref{prop:orderedRandomCanonical}. Since every infinite orbit $X$ of  $(\P^+,c_1,\ldots,c_n)$ induces $\P^+$, $h$ must have one of the behaviours of Lemma~\ref{lem:behaviours} on $X$. By Lemma~\ref{lem:chainantichain}, we may assume that $h$ behaves like $\rev$ or like $\id$ on every infinite orbit of $(\P^+,c_1,\ldots,c_n)$; for otherwise, $\cl{\{f \}\cup\Aut(\P)}$ is the full symmetric group $\Sym_P$, which implies that $\{f,f\inv \}\cup\Aut(\P)$ M-generates all injective functions, and in particular a function with the desired properties.
	  
	   Now let $S\subseteq P$ be a $\prec$-clean skeleton of $(\P,c_1,\ldots,c_n)$.	 %Since $\{f, f\inv\}\cup\Aut(\P)$ does not generate all permutations on $P$, Lemma~\ref
%{lemma:ruleOutChainAntichain} implies that $h$ behaves like $\id$ or like $\rev$ on all infinite orbits of $
%(\P,c_1,\ldots,c_n)$. 
We claim that $h$, considered as a function from $(\P,c_1,\ldots,c_n)$ to $(\P,d_1,\ldots,d_m)$,  is canonical on $S$, that is, it satisfies the definition of canonicity for tuples in $S$. To see this, let $a=(a_1,a_2),b=(b_1,b_2)\in S^2$ have the same type in  $(\P,c_1,\ldots,c_n)$. Then either $a,b$ or $a,\tilde b=(b_2,b_1)$ have the same type in $(\P^+,c_1,\ldots,c_n)$, and so either $h(a),h(b)$ or $h(a),h(\tilde b)$ have the same type in $(\P^+,d_1,\ldots,d_m)$, and hence also in $(\P,d_1,\ldots,d_m)$. In the first case we are done; in the second case, $\tp(a)=\tp(b)=\tp(\tilde b)$ in  $(\P,c_1,\ldots,c_n)$ implies that $a_1,a_2,b_1,b_2$ all belong to the same orbit in $(\P,c_1,\ldots,c_n)$. Since $h$ behaves like $\rev$ or like $\id$ on this orbit, we conclude that $f(a),f(b)$ have the same type in $(\P,d_1,\ldots,d_m)$.

Let $i:(P;\leq,c_1,\ldots,c_n)\To (S;\leq,c_1,\ldots,c_n)$ be an isomorphism, and set $g:=h\circ i$. Then $g$ is canonical as a function from $(\P,c_1,\ldots,c_n)$ to $(\P,d_1,\ldots,d_m)$, and agrees with $f$ on $\{c_1,\ldots,c_n\}$. Since $i$ preserves $\leq$ and its negation, it is M-generated by $\Aut(\P)$. Hence so is $g$, proving the lemma.
\end{proof}
\section{Applying canonical functions}

\subsection{Ordering orbits}

\begin{definition}
For disjoint subsets $X,Y$ of $P$ we write 
\begin{itemize}
	\item $X\leq Y$ iff there exist $x\in X$, $y\in Y$ such that $x\leq y$;
	\item $X\bot Y$ iff $x\bot y$ for all $x\in X$, $y\in Y$;
	\item $X<Y$ iff $x<y$ for all $x\in X$ and all $y\in Y$.
\end{itemize}
We call $X,Y$ \emph{incomparable} iff $X\bot Y$, and \emph{comparable} otherwise (which is the case iff $X\leq Y$ or $Y\leq X$). We say that $X,Y$ are \emph{strictly comparable} iff $X<Y$ or $Y<X$.
\end{definition}

\begin{lem}
Let $c_1,\ldots,c_n\in P$. The relation $\leq$ defines a partial order on the orbits of $(\P,c_1,\ldots,c_n)$.
\end{lem}
\begin{proof}
Reflexivity is obvious. To see that $X\leq Y$ and $Y\leq X$ imply $X=Y$, observe first that it follows from Fact~\ref{fact:orbits} that $X$ is convex, i.e., if $x,z\in X$ satisfy $x\leq z$ and $y\in P$ is so that 
$x\leq y$ and $y\leq z$, then $y\in X$. Now there exist $x, x'\in X$ and $y,y'\in Y$ such that $x\leq y$ and $x'\geq y'$. Since $y,y'$ belong to the same orbit, they satisfy the same first-order formulas over $(\P,c_1,\ldots,c_n)$, and hence there exists $z\in X$ such that $z\geq y$. Since $X$ is convex, we have $y\in X$, which is only possible if $X=Y$ since distinct orbits are disjoint.

Suppose that $X\leq Y$ and $Y\leq Z$. Then there exist $x\in X$, $y,y'\in Y$ and $ z\in Z$ such that $x\leq y$ and $y'\leq z$. Since $y,y'$ satisfy the same first-order formulas, there exists $x'\in X$ such that $x'\leq y'$. Hence $x'\leq z$ and so $X\leq Z$, proving transitivity.
\end{proof}

Let $X, Y$ be infinite orbits of $(\P,c_1,\ldots,c_n)$. Then precisely one of the following cases holds.
\begin{itemize}
	\item $X$ and $Y$ are strictly comparable;
	\item $X$ and $Y$ are incomparable;
	\item $X$ and $Y$ are comparable, but not strictly comparable.
\end{itemize}

In the third case, if $X\leq Y$, then there exist $x,x'\in X$ and $y,y'\in Y$ such that $x<y$ and $x'\bot y'$, and there are no $x''\in X$ and $y''\in Y$ such that $x''>y''$.

\begin{definition}
	If for two disjoint subsets $X,Y$ of $P$ we have $X\leq Y$, $Y\nleq X$, and $X\not < Y$, or vice-versa, then we write $X\div Y$.
\end{definition}

\subsection{Behaviors generating $\Sym_P$}

\begin{definition}
	Let $X,Y\subseteq P$ be disjoint, and let $f:P\To P$ be a function. We say that $f$
	\begin{itemize}
		\item \emph{behaves like $\id$ on $X$} iff $x< x'$ implies $f(x)<f(x')$ and $x\bot x'$ implies $f(x)\bot f(x')$ for all $x,x'\in X$;
		\item \emph{behaves like $\rev$ on $X$} iff $x< x'$ implies $f(x)>f(x')$ and $x\bot x'$ implies $f(x)\bot f(x')$ for all $x,x'\in X$;
		\item
		\emph{behaves like $\id$ between $X$ and $Y$} iff $x< y$ implies $f(x)<f(y)$, $x>y$ implies $f(x)>f(y)$, and $x\bot y$ implies $f(x)\bot f(y)$ for all $x\in X, y\in Y$.
		\end{itemize}
\end{definition}

\begin{lem}\label{lem:idOrrev}
Let $\G\supseteq \Aut(\P)$ be a closed group, and let $c_1,\ldots,c_n\in P$. Let $g:(\P,c_1,\ldots,c_n)\To\P$ be a canonical function M-generated by $\G$. Then $g$ behaves like $\id$ or like $\rev$ on each infinite orbit $X$ of $(\P,c_1,\ldots,c_n)$, or else $\G=\Sym_P$. 
\end{lem}
\begin{proof}
	Let $X$ be an infinite orbit, and let $x,x'\in X$ such that $x\bot x'$. Then the type of $(x,x')$ in $(\P,c_1,\ldots,c_n)$ equals the type of $(x',x)$ in $(\P,c_1,\ldots,c_n)$. Hence, the type of $(g(x),g(x'))$ must equal the type of $(g(x'),g(x))$ in $\P$, which is only possible if $g(x)\bot g(x')$, and hence $g$ preserves $\bot$ on $X$.
	
	Now if $g(a)<g(a')$ for some $a,a'\in X$ with $a<a'$, then the same holds for all $a,a'\in X$ with $a<a'$, and $g$ behaves like $\id$ on $X$. If $g(a')<g(a)$ for some $a,a'\in X$ with $a<a'$, then $g$ behaves like $\rev$ on $X$. Finally, if $g(a)\bot g(a')$ for some $a,a'\in X$ with $a<a'$, then $g$ sends $X$ to an antichain. Since $X$ contains all finite partial orders, and by the homogeneity of $\P$, we can then refer to Lemma~\ref{lem:chainantichain} to conclude that $\G=\Sym_P$.
\end{proof}

\begin{lem}\label{lem:preservesDiv}
Let $\G\supseteq \Aut(\P)$ be a closed group, and let $c_1,\ldots,c_n\in P$. Let $g:(\P,c_1,\ldots,c_n)\To\P$ be a canonical function M-generated by $\G$. Then $g[X]\div g[Y]$ for all infinite orbits $X,Y$ of $(\P,c_1,\ldots,c_n)$ with $X\div Y$, or else $\G=\Sym_P$. 
\end{lem}
\begin{proof}
	Suppose there are infinite orbits $X,Y$ with $X\div Y$ but for which $g[X]\div g[Y]$ does not hold. Assume without loss of generality that $X\leq Y$. By Lemma~\ref{lem:idOrrev}, we may assume that $g$ behaves like $\id$ or like $\rev$ on $X$ and on $Y$.
	
	First consider the case where $g[X]<g[Y]$ or $g[Y]<g[X]$. Let $A\subseteq P$ be finite; we claim that $\G$ M-generates a function which sends $A$ to a chain. There is nothing to show if $A$ is itself a chain, so assume that there exist $x,y$ in $A$ with $x\perp y$. Then using the extension property, one readily checks that there exists $\alpha\in\Aut(\P)$ which sends the principal ideal of $x$ in $A$ into $X$ and all other elements of $A$, and in particular $y$, into $Y$. Set $h:=g\circ \alpha$. Then $h(x)$ and $h(y)$ are comparable, and $h$ does not add any incomparabilities between elements of $A$. Hence, repeating this procedure and composing the functions, we obtain a function which sends $A$ to a chain. Lemma~\ref{lem:chainantichain} then implies $\G=\Sym_P$. 

The other case is where $g[X]\bot g[Y]$. Then an isomorphic argument shows that we can map any finite subset $A$ of $P$ to an antichain via a function which is M-generated by $\G$. Again, Lemma~\ref{lem:chainantichain} yields $\G=\Sym_P$.
\end{proof}

\begin{lem}\label{lem:sameOnAll}
Let $\G\supseteq \Aut(\P)$ be a closed group, and let $c_1,\ldots,c_n\in P$. Let $g:(\P,c_1,\ldots,c_n)\To\P$ be a canonical function M-generated by $\G$. Then $g$ behaves like $\id$ on all infinite orbits of $(\P,c_1,\ldots,c_n)$, or it behaves like $\rev$ on all infinite orbits of $(\P,c_1,\ldots,c_n)$, or else $\G=\Sym_P$. 
\end{lem}
\begin{proof}
By Lemma~\ref{lem:idOrrev}, we may assume that $g$ behaves like $\id$ or $\rev$ on all infinite orbits.
Suppose that the behaviour of $g$ is not the same on all infinite orbits. Consider the graph $H$ on the infinite orbits of $(\P,c_1,\ldots,c_n)$ in which two orbits $X,Y$ are adjacent if and only if $X\div Y$ holds. We claim that $H$ is connected. To see this, let $X,Y$ be infinite orbits with $X<Y$. Pick $x,x'\in X$ and $y,y'\in Y$ such that $x<x'$ and $y'< y$. By the extension property, there exists $z\in P$ such that $x<z$, $z\bot x'$, $z\bot y'$, and $z<y$. Let $Z$ be the orbit of $z$ in $(\P,c_1,\ldots,c_n)$. Then $X\div Z$ and $Z\div Y$, and so there is a path from $X$ to $Y$ in $H$.  Now if $X,Y$ are infinite orbits which are incomparable, then there exists an infinite orbit $Z$ with $X<Z$ and $Y<Z$, and so again there is a path from $X$ to $Y$ in $H$.

Since $H$ is connected, there exist infinite orbits $X,Y$ with $X\div Y$ such that $g$ behaves like $\id$ on $X$ and like $\rev$ on $Y$. Assume that $X\leq Y$; the proof of the case $Y\leq X$ is dual. By Lemma~\ref{lem:preservesDiv}, we may furthermore assume that $g[X]\div g[Y]$, or else we are done. This leaves us with two possibilities, $g[X]\leq g[Y]$ or $g[Y]\leq g[X]$. 

The first case  $g[X]\leq g[Y]$ splits into two subcases:

\begin{itemize}
\item For all $x\in X$, $y\in Y$, $x<y$ implies $g(x)<g(y)$ and $x\bot y$ implies $g(x)\bot g(y)$;
\item For all $x\in X$, $y\in Y$, $x<y$ implies $g(x)\bot g(y)$ and $x\bot y$ implies $g(x)< g(y)$.
\end{itemize}

Let $x,x'\in X$ and $y,y'\in Y$ be so that $x< x'$, $x<y'$, $x'<y$, $y'<y$, and $x'\bot y'$. Then in the first subcase we can derive $g(x')<g(y)$, $g(y)<g(y')$, and $g(x')\bot g(y')$, a contradiction. In the second subcase, $g(x)<g(x')$, $g(x')<g(y')$,  and $g(x)\bot g(y')$, again a contradiction.

In the second case $g[Y]\geq g[X]$ we have the following possibilities:

\begin{itemize}
\item For all $x\in X$, $y\in Y$, $x<y$ implies $g(x)>g(y)$ and $x\bot y$ implies $g(x)\bot g(y)$;
\item For all $x\in X$, $y\in Y$, $x<y$ implies $g(x)\bot g(y)$ and $x\bot y$ implies $g(x)> g(y)$.
\end{itemize}

Let $x,x'\in X$ and $y,y'\in Y$ be as before. Then in the first subcase we can derive $g(x)<g(x')$, $g(y')<g(x)$, and $g(x')\bot g(y')$, a contradiction. In the second subcase, $g(y)<g(y')$, $g(y')<g(x')$,  and $g(y)\bot g(x')$, again a contradiction.
\end{proof}

\subsection{Behaviors generating $\R$}

\begin{lem}\label{lem:generatingRev}
Let $\G\supseteq \Aut(\P)$ be a closed group, and let $c_1,\ldots,c_n\in P$. Let $g:(\P,c_1,\ldots,c_n)\To\P$ be a canonical function M-generated by $\G$. If $g$ behaves like $\rev$ on some infinite orbit of $(\P,c_1,\ldots,c_n)$, then $\G\supseteq \R$.
\end{lem}
\begin{proof}
Let $X$ be the infinite orbit. Pick an isomorphism $i:(P;\leq)\to (X;\leq)$. Then given any finite $A\subseteq P$, there exists $\alpha\in\Aut(\P)$ such that $\alpha\circ g\circ i$ agrees with $\rev$ on $A$. Since $g$ and $i$ are generated by $\G$, there exists $\beta\in\G$ such that $\beta$ agrees with $\rev$ on $A$. Hence, $\rev\in\G$.
\end{proof}

\subsection{Behaviors generating $\T$}

\begin{lem}\label{lem:generatingTurn}
Let $\G\supseteq \Aut(\P)$ be a closed group, and let $c_1,\ldots,c_n\in P$. Let $g:(\P,c_1,\ldots,c_n)\To\P$ be a canonical function M-generated by $\G$ which behaves like $\id$ on all of its orbits. Then $g$ behaves like $\id$ between all infinite orbits of $(\P,c_1,\ldots,c_n)$, or else $\G\supseteq \T$. 
\end{lem}
\begin{proof}
	Let infinite orbits $X,Y$ be given. 	
		
	We start with the case $X\div Y$. Say without loss of generality $X\leq Y$. By Lemma~\ref{lem:preservesDiv}, we may assume that $g[X]\div g[Y]$, or else $\G=\Sym_P$. Hence $g[X]\leq g[Y]$ or $g[Y]\leq g[X]$. If $g[X]\leq g[Y]$, then either $g$ behaves like $\id$ between $X$ and $Y$ and we are done, or $x<y \rightarrow g(x)\bot g(y)$ and $x\bot y \rightarrow g(x)< g(y)$ hold for all $x\in X$, $y\in Y$; the latter, however, is impossible, as for $x,x'\in X$ and $y\in Y$ with $x< x'$, $x< y$, and $x'\bot y$ we would have $g(x)<g(x')<g(y)$ and $g(x)\bot g(y)$. Now suppose $g[Y]\leq g[X]$. Then we have one of the following:
	\begin{itemize}
\item For all $x\in X$, $y\in Y$, $x<y$ implies $g(x)>g(y)$ and $x\bot y$ implies $g(x)\bot g(y)$;
\item For all $x\in X$, $y\in Y$, $x<y$ implies $g(x)\bot g(y)$ and $x\bot y$ implies $g(x)> g(y)$.
	\end{itemize}
The first case is absurd since picking $x,x',y$ as above yields $g(x)<g(x')$, $g(x)>g(y)$, and $g(x')\bot g(y)$.  We claim that in the second case $\G$ contains $\turn$. Let $F\subseteq P$ be any random filter.  Let $A\subseteq P$ be finite, and set $A_2:=A\cap F$, and $A_1:=A\setminus A_2$. Then there exists an automorphism $\alpha$ of $\P$ which sends $A_2$ into $Y$ and $A_1$ into $X$. The composite $g\circ \alpha$ behaves like $\turn_F$ on $A$ for what concerns comparabilities and incomparabilities, and hence there exists $\beta\in\Aut(\P)$ such that $\beta\circ g\circ \alpha$ agrees with $\turn_F$ on $A$. By topological closure we infer $\turn_F\in \G$.
	
	Now consider the case where $X,Y$ are strictly comparable, say $X<Y$. Then we know from the proof of Lemma~\ref{lem:sameOnAll} that  there exists an infinite orbit $Z$ such that $X\leq Z\leq Y$,  $X\div Z$ and $Z\div Y$. Let $x\in X$ and $y\in Y$ be arbitrary. There exists $z\in Z$ such that $x<z<y$. As $g$ behaves like $\id$ between $X$ and $Z$ and between $Z$ and $Y$, we have that $g(x)<g(z)<g(y)$, and hence $g$ behaves like $\id$ between $X$ and $Y$.
	
	It remains to discuss the case $X\bot Y$. Suppose that $g[X]$ and $g[Y]$ are comparable, say $g[X]<g[Y]$. Then given any finite $A\subseteq P$ with incomparable elements $x,y$, using the extension property we can find $\alpha\in\Aut(\P)$ which sends $x$ into $X$, all elements of $A$ which are incomparable with $x$ into $Y$, and all other elements of $A$ into infinite orbits which are comparable with both $X$ and $Y$. Applying $g\circ \alpha$ then increases the number of comparabilities on $A$, and hence repeated applications of such functions will send $A$ onto a chain, proving $\G= \Sym_P$.
\end{proof}

\begin{lem}\label{lem:generatingTurn2}
Let $\G \supseteq \Aut(\P)$ be a closed group, and let $c_1,\ldots,c_n\in P$. Let $g:(\P,c_1,\ldots,c_n)\To\P$ be a canonical function M-generated by $\G$ which behaves like $\id$ on all of its orbits. Then $g$ behaves like $\id$ between all orbits of $(\P,c_1,\ldots,c_n)$ (including the finite ones), and hence is M-generated by $\Aut(\P)$, or else $\G\supseteq\T$. 
\end{lem}
\begin{proof}
Let $1\leq i\leq n$, and let $X$ be an infinite orbit which is incomparable with $\{c_i\}$. Suppose that $g[X]$ and $\{g(c_i)\}$ are strictly comparable, say $\{g(c_i)\}< g[X]$. Let $Y$ be an infinite orbit such that $X\leq Y$, $X\div Y$, and $\{c_i\}< Y$. Let moreover $Z$ be an infinite orbit such that $Z<\{c_i\}$, $Z\leq X$ and $Z\div X$. Then by the preceding lemma, we may assume that $g$ behaves like $\id$ between $X, Y$ and $Z$. We cannot have $g[Z]<\{g(c_i)\}$ as this would imply $g[Z]<g[X]$, contradicting the fact that $g$ behaves like $\id$ between $Z$ and $X$. Suppose that $g[Z]\bot \{g(c_i)\}$. Set $S:=Z\cup X\cup Y\cup\{c_i\}$. Then it is easy to see that $(S;\leq)$ satisfies the extension property, and hence is isomorphic which $\P$; fix an isomorphism $i:(P;\leq,c_i)\To (S;\leq,c_i)$. This isomorphism is M-generated by $\Aut(\P)$ since it can be approximated by automorphisms of $\P$ on all finite subsets of $P$. The restriction of $g$ to $S$ is canonical as a function from $(S;\leq,c_i)$ to $\P$. Hence, the function $h:= g\circ i$ is canonical as a function from $(\P,c_i)$ to $\P$, and has the same behaviour as the restriction of $g$ to $S$. Let $\alpha\in\Aut(\P)$ be so that $\alpha(h(c_i))=c_i$. Then $t:=h\circ\alpha\circ h$ has the property that $t(x)>t(c_i)$ for all $x\neq c_i$, and that $t(x)\bot t(y)$ if and only if $x\bot y$, for all $x,y\in P\setminus\{c_i\}$. Hence, given any finite $A\subseteq P$ which is not a chain, we can pick $x\in A$ which is not comparable to all other elements of $A$, and find $\beta\in\Aut(\P)$ which sends $x$ to $c_i$; then $t\circ\beta$ strictly increases the number of comparabilities among the elements of $A$. Repeating this process and composing the functions, we find a function which is M-generated by $\G$ and which maps $A$ onto a chain. Hence, $\G=\Sym_P$.

Therefore, we may henceforth assume that $g$ behaves like $\id$ between all $\{c_i\}$ and all infinite orbits $X$ with $\{c_i\}\bot X$. Now suppose that there exists $1\leq i\leq n$ and an infinite orbit $X$ with $X<\{c_i\}$ such that $\{g(c_i)\}<g[X]$. Pick an infinite orbit $Y$ which is incomparable with $c_i$, and which satisfies $X\leq Y$. Then $\{g(c_i)\}<g[Y]$ since $g$ behaves like $\id$ between $X$ and $Y$, a contradiction. Next suppose there exists $1\leq i\leq n$ and an infinite orbit $X$ with $X<\{c_i\}$ such that $\{g(c_i)\}\bot g[X]$. Then pick an infinite orbit $Y$ as in the preceding case, and an infinite orbit $Z$ with $\{c_i\}<Z$. Now given any finite $A\subseteq P$ which does not induce an antichain, we can pick $y\in A$ which is not minimal in $A$. Taking $\alpha\in\Aut(\P)$ which sends $y$ to $c_i$ and $A$ into $X\cup Y\cup Z\cup\{c_i\}$, we then have that application of $g\circ\alpha$ increases the number of incomparabilites of $A$. Repeated composition of such functions yields a function which sends $A$ onto an antichain. Hence, $\G=\Sym_P$. The case where there  exist $1\leq i\leq n$ and an infinite orbit $X$ with $\{c_i\}<X$ such that $\{g(c_i)\}\bot g[X]$ is dual.

We turn to the case where we have two distinct finite orbits $\{c_i\}$ and $\{c_j\}$. Suppose first that they are comparable, say $c_i<c_j$. Picking an infinite orbit $Z$ with $\{c_i\}< Z<\{c_j\}$ then yields, by what we know already, $\{g(c_i)\}<g[Z]<\{g(c_j)\}$, so we are done. Finally, suppose that $c_i\bot c_j$. Then given any finite $A\subseteq P$ which has incomparable elements $x,y$, we can send $x$ to $c_i$, $y$ to $c_j$, and the rest of $A$ to infinite orbits via some $\alpha\in\Aut(\P)$. But then application of $g\circ \alpha$ increases the number comparabilities on $A$, and hence repeating the process yields a function which sends $A$ to a chain. Hence, $\G=\Sym_P$.
\end{proof}

\subsection{Climbing up the group lattice}

\begin{proposition}\label{prop:aboveP}
Let $\G\supsetneq \Aut(\P)$ be a closed group. Then $\G$ contains either $\R$ or $\T$.
\end{proposition}
\begin{proof}
There exist $\pi\in \G\setminus\Aut(\P)$ and elements $u,v\in P$ such that $u\leq v$ and $\pi(u) \nleq \pi(v)$. Let $g:(\P,u,v)\To \P$ be a canonical function M-generated by $\G$ which agrees with $\pi$ on $\{u,v\}$. If $g$ behaves like $\rev$ on some infinite orbit of $\Puv$, then $\G\supseteq\R$ by Lemma~\ref{lem:generatingRev}. Otherwise Lemma~\ref{lem:generatingTurn2} states that $g$ is generated by $\Aut(\P)$ or $\G\supseteq\T$. Since $g(u) \nleq g(v)$, only the latter possibility can be the case.
\end{proof}

\begin{proposition}\label{prop:aboveR}
Let $\G\supsetneq \R$ be a closed group. Then $\G$ contains $\T$.
\end{proposition}
\begin{proof}
Let $\pi\in \G\setminus \R$. Then there exists a finite tuple $c=(c_1,\ldots,c_n)$ of elements of $P$ such that no function in $\R$ agrees with $\pi$ on $c$.  Let $g:(\P,c_1,\ldots,c_n)\to \P$ be a canonical function which is M-generated by $\G$ and which agrees with $\pi$ on $\{c_1,\ldots,c_n\}$. By Lemma~\ref{lem:sameOnAll}, we may assume that either $g$ behaves like $\id$ on all infinite orbits, or it behaves like $\rev$ on all infinite orbits of $(\P,c_1,\ldots,c_n)$. By composing $g$ with $\rev$, we may assume that it behaves like $\id$ on all infinite orbits. But then Lemma~\ref{lem:generatingTurn2} implies that $\G\supseteq\T$, or that $g$ is M-generated by $\Aut(\P)$. The latter is, of course, impossible.
\end{proof}

\subsection{Relational descriptions of $\T$ and $\M$} Before climbing up further, we need to describe the groups $\T$ and $\M$ relationally. 
The componentwise action of the group $\T$ on triples of distinct elements of $P$ has three orbits, namely:

\begin{itemize}
\item[$\oo$:] the orbit of  the 3-element
  antichain, i.e., the set of all tuples $(a,b,c)\in P^3$ such that one of the following holds: $a\perp b, b\perp c, c\perp a$;

$a<b,a<c, b  \perp c $; \quad $b<a,b<c, a  \perp c $; \quad $c<a,c<b, b  \perp c $;

$a>b,a>c, b  \perp c $; \quad $b>a,b>c, a  \perp c $; \quad $c>a,c>b, b  \perp c $; 

\item[$\oz$:] the orbit of the 3-element chain $a<b<c$, i.e., the set of all $(a,b,c)\in P^3$ such that one of the following holds: 

$a<b<c$;\quad $b<c<a$;\quad $c<a<b$;

$a <b, c\perp a, c\perp b$;\quad
$b <c, a\perp b, a\perp c$;\quad
$c <a, b\perp a, b\perp c$;
\item[$\od$:] the dual of $\oz$; that is, the orbit of the chain $a>b>c$, or more precisely the set of all $(a,b,c)\in P^3$ such that one of the following holds: 

$a>b>c$; \quad $b>c>a$; \quad $c>a>b$;

$a >b, c\perp a, c\perp b$;\quad
$b >c, a\perp b, a\perp c$;\quad
$c >a, b\perp a, b\perp c$.
\end{itemize}

\begin{defn}\label{defn:rotations}
Let $\{X,Y,Z\}$ be a partition of $P$ into disjoint subsets such that $X$ is an ideal of $\P$, $Z$ is a filter of $\P$, $X\leq Y$, $Y \leq Z$ and $X<Z$. A \emph{rotation} on $\P$ with respect to $X,Y,Z$ is any permutation $f$ on $P$ which behaves like $\id$ on each class of the partition, and such that for all $x\in X$, $y\in Y$, and $z\in Z$ we have
\begin{itemize}
\item $f(z)<f(x)$;
\item $f(y)<f(x)$ iff $x\bot y$ and $f(y)\bot f(x)$ iff $x<y$;
\item $f(z)<f(y)$ iff $y\bot z$ and $f(z)\bot f(y)$ iff $y<z$.
\end{itemize}
\end{defn}

Observe that if $F$ is a random filter, then $\turn_F$ is a rotation with respect to the partition $\{\emptyset,P\setminus F,F\}$.

\begin{proposition}\label{prop:rotations}
$\T$ contains all rotations on $\P$.
\end{proposition}
\begin{proof}
Let $f$ be a rotation on $\P$, let $\{X,Y,Z\}$ be the corresponding partition, and let $S\subseteq P$ be finite. Set $X':= X\cap S$, $Y':= Y\cap S$, and $Z':=Z\cap S$. Let $F\subseteq P$ be a random filter with $F\supseteq Z'$ and $P\setminus F \supseteq X'\cup Y'$. Since $\turn_F(u)\not < \turn_F(z)$ for all $u\in X'\cup Y'$ and all $z\in Z'$, there exists a random filter $F'$ with $F'\supseteq \turn_F[X'\cup Y']$ and $P\setminus F' \supseteq \turn_F[Z']$. It is a straightforward verification that $\turn_{F'}\circ \turn_F$ changes the relations between elements of $X'\cup Y'\cup Z'$ in the very same way as the rotation $f$, and hence there exists an automorphism $\alpha$ of $\P$ such that $\alpha\circ\turn_{F'}\circ \turn_F$ agrees with $f$ on $X'\cup Y'\cup Z'$.
\end{proof}

\begin{lem}\label{lem:descr:turn}
	$\T=\Aut(P;\oo,\oz,\od)$.
\end{lem}
\begin{proof}
To show that $\turn$ preserves $\oo$, $\oz$ and $\od$ is only a matter of verification of a finite number of cases. For the converse, let $f\in \Aut(P;\oo,\oz,\od)$; we show it is a rotation. Define a binary relation $\sim$ on $P$ by setting $x\sim y$ if and only if $(x,y)$ and $(f(x),f(y))$ have the same type in $\P$, for all $x,y\in P$. Clearly, $\sim$ is reflexive and symmetric; we claim it is transitive, and hence an equivalence relation. To this end, let $x,y,z\in P$ such that $x\sim y$ and $y\sim z$. Now by going through all possible relations that might hold between $x,y,z$, using the fact that these relations remain unaltered between $x$ and $y$ as well as between $y$ and $z$, and taking into account the fact that $(x,y,z)$ in $\oo$ ($\oz$, $\od$) implies $(f(x),f(y),f(z))$ in $\oo$ ($\oz$, $\od$), one checks that the relation which holds between $x$ and $z$ has to remain unchanged as well -- this is a finite case analysis which we leave to the reader.

If $\sim$ has only one equivalence class, then $f$ it is an automorphism of $\P$ and there is nothing to show, so assume henceforth that this is not the case. Then there exist equivalence classes $X$, $Y$ and $x\in X$, $y\in Y$ such that $x\bot y$; we may assume without loss of generality that $f(x)> f(y)$. 

Let $u,v\in X\cup Y$ such that $u<v$, and suppose that $f(v)<f(u)$. Pick $r\in P$ incomparable with $u,v,x,y$. Then $(r,x,y)\in\oo$, so $(f(r),f(x),f(y))\in\oo$. Consequently, $f(r)\bot f(x)$ or $f(r)\bot f(y)$, and hence $r\in X\cup Y$. 
 Now observe that $(u,v,r)$, and hence also its image under $f$, is an element of $\oz$. 
Hence $f(v)<f(u)$ yields $f(v)<f(r)<f(u)$, contradicting $r\in X\cup Y$. We conclude that comparable elements of $X\cup Y$ either belong to the same class, or they are sent to incomparable elements.

Pick any $u\in P$ such that $u<x$ and $u\bot y$. Then $(f(u),f(x),f(y))\in\oz$ and $f(x)>f(y)$ imply $f(u)<f(x)$, and so $u\in X$. Similarly, any $v\in P$ such that $y<v$ and $v \bot x$ is an element of $Y$, and in particular $X\leq Y$.

We next claim that $Y\nleq X$. Suppose there exist $u\in Y$, $v\in X$ with $u<v$. If $u>x$, then $(x,u,v)\in\oz$, and so $f(x)<f(v)$  and the fact that we cannot have $f(x)<f(u)$ yield a contradiction. Hence, $u\not > x$, and by symmetry $v\not <y$.  Suppose $v>y$. If $v>x$, then $(y,x,v)\in\oo$, but $f(y)<f(x)<f(v)$, a contradiction. By the preceding paragraph, $v\bot x$ would imply $v\in Y$; so $v\bot y$, and by symmetry $u\bot x$. If $u<y$ and $x<v$, then $f(u)<f(y)<f(x)<f(v)$, contradicting the fact that $u$ and $v$ are elements of different classes. So assume without loss of generality that $u\not < y$; since $u> y$ would imply $v> y$, which we already excluded, we then have  $u\bot y$. Since $(x,u,y)\in\oo$ and $f(x)>f(y)$, we conclude $f(u)<f(x)$. Hence, if $v>x$, then $f(u)<f(x)<f(v)$, a contradiction, so we must have $v\bot x$. But then $(u,v,x)\in\oz$, $f(u)<f(x)$, and the fact that $f(v)$ is incomparable with $f(u)$ and $f(x)$ yield the final contradiction.

Suppose there exist $u\in X$ and $y\in Y$ with $u\bot v$ and such that $f(u)<f(v)$. As above, we could then conclude that $X\nleq Y$, a contradiction.

Say that $A,B$ are equivalence classes for which $A<B$. Picking $a\in A$, $b\in B$, and any $c\in P$ which is incomparable with $a$ and $b$, we then have $(a,b,c)\in\oz$. We cannot have $c\in A\cup B$, and so $f(c)$ must be comparable with $f(a)$ and $f(b)$. The only possibility then is that $f(b)<f(a)$.

Let $Z$ be an equivalence class distinct from $X,Y$ and such that $Y\leq Z$. Then $X\leq Z$. We claim that $Z>Y$ is impossible. Otherwise, there exist $x\in X$, $y\in Y$, and $z\in Z$ such that $x<y<z$, and so $(x,y,z)\in \oz$. But $f(x)\bot f(y)$ and $f(z)<f(y)$ imply $(f(x),f(y),f(z))\nin \oz$, a contradiction. We next claim that $X<Z$. Otherwise, pick $x\in X$ and $z\in Z$ with $x\bot z$, and an arbitrary $y\in Y$ such that $x< y$. 
Then $(f(x),f(y),f(z))\in\oz$, $f(x)> f(z)$ and $f(x)\bot f(y)$ yield a contradiction. Suppose next that there exist two distinct classes $Z_1,Z_2$ with $Y\leq Z_1,Z_2$. We know that $Z_1,Z_2$ must be comparable, say $Z_1\leq Z_2$. Pick $z_1\in Z_1$, $z_2\in Z_2$ with $z_1<z_2$. Since $X< Z_1, Z_2$, we then have $f(x)>f(z_1),f(z_2)$, and $f(z_1)\not < f(z_2)$ yields a contradiction. So there is at most one class $Z$ distinct from $Y$ with $Z\geq Y$, and it satisfies $Z>X$ and $Z\not > Y$. 

Similarly there is at most one class $W$ distinct from $X$ with $W\leq X$, and it satisfies $W<Y$ and $W\not < X$. By the same kind of argument that yielded uniqueness of $Z$ above, 
$W$ and $Z$ cannot exist simultaneously, say that $W$ does not. Let $U$ be any other class distinct from $X$ and $Y$. 
Then $X\leq U\leq Y$, and so $X<Y$, a contradiction.

If $Z$ does not exist, then $Y$ is a filter and $f$ is of the form $\turn_Y$. 
If $Z$ does exist, then $f$ is a rotation with respect to the partition $\{X,Y,Z\}$.
\end{proof}

\begin{cor}\label{cor:turnIsRotations}
The group $\T$ consists precisely of the rotations on $\P$. In particular, the composition of two rotations is again a rotation.
\end{cor}
\begin{proof}
 By Lemma~\ref{lem:descr:turn}, if $f\in \T$, then $f\in \Aut(P;\oo,\oz,\od)$. It then follows from the proof of the other direction of same lemma that $f$ is a rotation.
\end{proof}

\begin{prop}\label{prop:descr:turn}
$\T=\Aut(P;\cyc)$.
\end{prop}
\begin{proof}
By Lemma~\ref{lem:descr:turn}, $\T\subseteq \Aut(P;\cyc)$. If the two groups were not equal, then $\Aut(P;\cyc)$ would contain a function $f$ which sends a triple $a=(a_1,a_2,a_3)$ in $\oo$ to a triple in $\od$. Moreover, by first applying a function in $\T$, we could assume that $a$ induces an antichain in $\P$. But then for any automorphism $\alpha$ of $\P$ sending $a$ to $(a_3,a_2,a_1)$ we would get that $f\circ\alpha$ sends $a$ to a triple in $\oz$, a contradiction.
\end{proof}

\begin{lem}\label{lem:switch}
	Let $f\in\Aut(P;\oo)\setminus\T$. Then for all $a\in P^3$ we have $a\in \oz$ if and only if $f(a)\in \od$, i.e., $f$ switches $\oz$ and $\od$.
\end{lem}
\begin{proof}
Suppose there exists $a=(a_1,a_2,a_3)\in\oz$ with  $f(a)\in \oz$ -- we will derive a contradiction, implying $f(a)\in\od$. By symmetry, it then follows that all tuples in $\od$ are sent to $\oz$, and we are done.

Since $f\in\Aut(P;\oo)\setminus\T$, there exists $b=(b_1,b_2,b_3)$ in $\oz$ such that  $f(b)\in\od$. We first claim that by replacing $a$ and $b$ with adequate triples, we may assume that both $a$ and $b$ are strictly ascending, i.e., $a_1<a_2<a_3$ and $b_1<b_2<b_3$. Otherwise, either all strictly ascending triples are sent to $\oz$, or all strictly ascending triples are sent to $\od$. Assume without loss of generality the former. Let $g\in\T$ be so that it sends some strictly ascending triple $e\in P^3$ to $b$. Then $f\circ g$ sends $e$ to $f(b)\in\od$; on the other hand, since $g$ is a rotation by Corollary~\ref{cor:turnIsRotations}, it sends some other strictly ascending triple $w\in P^3$ onto a strictly ascending triple, and so $f\circ g(w)\in\oz$. Thus by replacing $f$ by $f\circ g$, $a$ by $w$ and $b$ by $e$, we may indeed henceforth assume that both $a$ and $b$ are strictly ascending triples.

Now let $c=(c_1,c_2,c_3)$ be a strictly ascending triple such that $a_i<c_j$ and $b_i<c_j$ for all $1\leq i,j\leq 3$. If $f(c)\in\oz$, then we replace $a$ by $c$, and otherwise we replace $b$ by $c$. Assume without loss of generality the former; hence, from now on we assume $b_1<b_2<b_3<a_1<a_2<a_3$, $f(b)\in\od$, and $f(a)\in\oz$. By replacing $f$ by $h\circ f$ for an appropriate function $h\in \T$ we may moreover assume that $f(a_i)=a_i$ for all $1\leq i\leq 3$.

Suppose that $f(b_i)\bot a_j$ for some $1\leq i,j\leq 3$. Then, for any $1\leq k\leq 3$ with $k\neq j$, the fact that $(b_i,a_j,a_k)\notin \oo$ implies $(f(b_i),a_j,a_k)\notin \oo$, and consequently $f(b_i)\bot a_k$. Hence, if $f(b_i)$ is incomparable with some $a_j$, then it is incomparable with all $a_j$, and if it is comparable with some $a_j$, then it is comparable with all $a_j$. Suppose that $f(b_i)\bot a_1$ for some $1\leq i\leq 3$, and consider $f(b_j)$, where $j\neq i$. Since $(b_i,b_j,a_1)\notin \oo$, we have  $(f(b_i),f(b_j),a_1)\notin \oo$. This implies that if $f(b_j)\bot f(b_i)$, then $f(b_j)$ and $a_1$ are comparable. Putting this information together, we conclude that any two distinct elements $f(b_i),f(b_j)$ which are incomparable with the $a_k$ are mutually comparable. Thus, the image of $S:=\{a_1,a_2,a_2,b_1,b_2,b_3\}$ under $f$ is the disjoint union of at most two chains; by applying $\turn_F$ for an appropriate random filter $F\subseteq P$, we may assume its image is a single chain. By the same argument, we may assume that $a_3$ is the largest element of this chain.

Since $f(b)\notin\oz$, there exists $b_i,b_j$ with
$b_i<b_j$ such that $f(b_j)<f(b_i)$. As in the following, we will not make use of the third element of $b$ anymore, we may assume that this is the case for $b_1,b_2$. Then either $f(b_2)<f(b_1)<a_2<a_3$, or $a_1<a_2<f(b_2)<f(b_1)$, or $f(b_2)<a_2<f(b_1)<a_3$. We will derive a contradiction from each of the three cases.

Pick any $u_1,u_2,u_3,u_4\in P$ such that $u_1<u_2$, $u_3<u_4$, and such that any other two elements $u_i,u_j$ are incomparable. Then there is an random filter $F\subseteq P$ containing $u_1,u_2$ but not $u_3,u_4$, and so $\turn_F(u_1)<\turn_F(u_2)<\turn_F(u_3)<\turn_F(u_4)$. Now if $f(b_2)<f(b_1)<a_2<a_3$, then by applying an automorphism of $\P$, we may assume that $(b_1,b_2,a_2,a_3)$ coincides with the ascending 4-tuple $t$ containing the $\turn_F(u_i)$. Picking an random filter $F'\subseteq P$ containing $a_2$ but not $f(b_1)$ and setting  $h:=\turn_{F'}\circ f\circ \turn_F$, we get that $h(u_2)<h(u_1)$, $h(u_3)<h(u_4)$, and all other $h(u_i), h(u_j)$ are incomparable. Pick any $x\in P$ such that $x>u_1$, $x> u_3$, $x\bot u_3$, and $x\bot u_4$. Then $(u_4,u_3,x)\in \oo$ implies that $h(x)>h(u_3)$, $(x,u_1,u_2)\in\oo$ implies $h(x)<h(u_1)$, and hence $h(u_3)<h(u_1)$, a contradiction.  If $a_1<a_2<f(b_2)<f(b_1)$, then by applying an automorphism of $\P$, we may assume that $(b_1,b_2,a_1,a_2)$ coincides with the tuple $t$. Picking an random filter $F'\subseteq P$ containing $f(b_2)$ but not $a_2$ and setting  $h:=\turn_{F'}\circ f\circ \turn_F$, we get that $h(u_2)<h(u_1)$, $h(u_3)<h(u_4)$, and all other $h(u_i), h(u_j)$ are incomparable, leading to the same contradiction as in the preceding case. Finally, assume $f(b_2)<a_2<f(b_1)<a_3$, and assume that $(b_1,b_2,a_2,a_3)$ coincides with $t$. Picking an random filter $F'\subseteq P$ containing $f(b_1)$ but not $a_2$ and setting  $h:=\turn_{F'}\circ f\circ \turn_F$, we get that $h(u_2)<h(u_3)$, $h(u_1)<h(u_4)$, and all other $h(u_i), h(u_j)$ are incomparable. Now pick $x\in P$ such that $x> u_i$ for all $1\leq i\leq 4$. Then $(u_1,u_2,x)\notin \oo$ implies that $h(x)\bot h(u_1)$ or $h(x)\bot h(u_2)$. However, $(u_2,u_4,x)\in\oo$ implies that $h(x)$ is comparable with $h(u_2)$, and similarly $(u_1,u_3,x)\in\oo$ implies that $h(x)$ is comparable with $h(u_1)$, a contradiction.

\end{proof}

\begin{prop}\label{prop:descr:max}
	$\M=\Aut(P;\pari)$.
\end{prop}
\begin{proof}
By Lemma~\ref{lem:descr:turn}, $\T$ is contained in $\Aut(P;\oo)$. Obviously, $\rev$ preserves $\oo$, so that indeed $\M\subseteq \Aut(P;\oo)$.

For the other direction, let $f\in \Aut(P;\oo)$. If $f\in\T$ then $f\in\M$ by definition of $\M$, so assume $f\nin \T$. Then $f$ switches $\oz$ and $\od$ by Lemma~\ref{lem:switch}. Since $\rev$ switches $\oz$ and $\od$ as well, $\rev\circ f$ preserves $\oo$, $\oz$ and $\od$. Thus, by Lemma~\ref{lem:descr:turn}, $\rev\circ f$ is an element of $\T$, and so $f\in\M$.
\end{proof}

\subsection{Climbing to the top}

\begin{proposition}\label{prop:aboveM}
Let $\G\supsetneq \M$ be a closed group. Then $\G$ is 3-transitive.
\end{proposition}
\begin{proof}
Since $\G$ is not contained in $\M$, $\oo$ cannot be an orbit of its componentwise action on $P^3$. Since it contains $\M$, the orbits of this action are unions of the orbits of the corresponding action of $\M$. However, the latter action has only two orbits of triples of distinct elements, namely $\oo$ and $\oz\cup\od$. Hence, $\G$ has only one such orbit, and is 3-transitive.
\end{proof}

\begin{proposition}\label{prop:3transitive}
Let $\G$ be a $3$-transitive closed group containing $\T$. Then $\G=\Sym_P$.
\end{proposition}

\begin{proof}
We prove by induction that $\G$ is $n$-transitive for all $n\geq 3$. Our claim holds for $n=3$ by assumption. So let $n\geq 4$ and assume that $\G$ is $(n-1)$-transitive. We claim that every $n$-element subset of $P$ can be mapped onto an antichain by a permutation in $\G$; $n$-transitivity then follows as in the proof of Lemma~\ref{lem:chainantichain}. We prove this claim in several steps, and will need the following partial orders.

For every natural number $k$ with $1\leq k\leq n$, let 
\begin{itemize}
	\item $S_n^k$ be the $n$-element poset consisting of $k$ independent points and a chain of $(n-k)$ elements below them;
	\item $T_n^k$ be the dual of $S_n^k$;

	\item  $A_n^k$ be the $n$-element poset consisting of $k$ independent points, an element below them, and an antichain of size $(n-k-1)$ independent from these points;
	\item $B_n^k$ be the dual of $A_n^k$;
	\item $C_k$ be the $k+1$-element poset consisting of $k$ independent points and an element below them; that is, $C_k=A_{k+1}^k= S_{k+1}^k$.
\end{itemize}

\textbf{Step 1: From anything to $A_n^k$ or $B_n^k$ for $k\geq \frac{n-1}{2}$.}	
	
	We first show that any $n$-element set $A\subseteq P$ can me mapped to a copy of $A_n^k$ or $B_n^k$, where $k\geq \frac{n-1}{2}$, via a function in $\G$. Let $A$  
	 be given, and write $A=A'\cup\{a\}$, where $A'$ has $n-1$ elements. Then by the induction hypothesis there exists $\pi\in\G$ which maps $A'$ to an antichain. Let $F\subseteq P$ be an random filter which separates $\pi(a)$ from $\pi[A']$, i.e., for all $b\in\pi[A']$ we have $b\in F$ if and only if $\pi(a)\nin F$. Then one can check that either $\pi[A]$ or $(\turn_F\circ\pi)[A]$ induce $A_n^k$ or $B_n^k$ in $\P$ for some $k\geq \frac{n-1}{2}$.

\textbf{Step 2: From $A_n^k$ ($B_n^k$) to $S_n^k$ ($T_n^k$) for $k\geq \frac{n-1}{2}$.}	

We now show that any copy of $A_n^k$ in $\P$ can be mapped to a copy of $S_n^k$ via a function in $\G$. The dual proof then shows that any copy of $B_n^k$ can be mapped to a copy of $T_n^k$.

Let $\{x_1, \ldots, x_{n-1}\}$ and $\{y_1,\ldots, y_{n-1}\}$ be disjoint subsets of $P$ inducing an antichain and a chain, respectively. By the $(n-1)$-transitivity of $\G$, the map $x_i\mapsto y_i$, $1\leq i\leq n-1$, can be extended to a permutation $\pi\in \G$. Let $X$ be the orbit of $(\P,x_1, \ldots, x_{n-1})$ such that $x\bot x_i$ for all $x\in X$ and all $1\leq i\leq n-1$. By Lemma~\ref{lem:interface} there exists a canonical function $g: (\P, x_1, \ldots, x_{n-1})\rightarrow (\P, y_1, \ldots, y_{n-1})$ M-generated by $\G$ that agrees with $\pi$ on $\{x_1, \ldots, x_{n-1}\}$. We may assume that $g$ behaves like $\id$ or like $\rev$ on $X$, by Lemma~\ref{lem:idOrrev}. If $g$ behaves like $\rev$ on $X$, then $\G$ contains $\rev$ by Lemma~\ref{lem:generatingRev}; replacing $g$ by $\rev\circ g$ and replacing each $y_i$ by $\rev(y_i)$, we may assume that $g$ behaves like $\id$ on $X$. Let $D\subseteq X$ be so that it induces $C_k$, and observe that $D':=D\cup\{x_1,\ldots,x_{n-k-1}\}$ induces a copy of $A_n^k$ in $\P$.
Since $g$ is canonical, all elements of $X$, and in particular all elements of $D$ are sent to the same orbit $Y$ of $(\P, y_1, \ldots, y_{n-1})$. Thus for all $1\leq i\leq n-1$ we have that either $g[D]<\{y_i\}$, or $g[D]\bot\{y_i\}$, or $g[D]>\{y_i\}$. Let $S$ be the set of those $y_i$ for which the first relation holds, and set $E:=g[D]\cup(\{y_1,\ldots,y_{n-1}\}\setminus S)$. Let $F\subseteq P$ be an random filter which separates $E$ from $S$, i.e., $F$ contains $S$, but does not intersect $E$. Then $\turn_F[S]\bot \turn_F[E]$. Choose an random filter $F'$ which contains  $\turn_F[S]$ and which does not intersect $\turn_F[E]$. Then $\turn_{F'}\circ\turn_F[S]< \turn_{F'}\circ\turn_F[E]$. Set $h:=\turn_F'\circ\turn_F\circ g$. Now for all $1\leq i\leq n-1$ we have that either $h[D]> \{h(x_i)\}$ or $h[D]\bot\{h(x_i)\}$. Moreover, $h$ behaves like $\id$ on $D$, and the $h(x_i)$ form a chain. Either there are at least $\frac{n-1}{2}$ elements among the $h(x_i)$ for which $h[D]> \{h(x_i)\}$, or there are at least $\frac{n-1}{2}$ of the $h(x_i)$ for which $h[D]\bot \{h(x_i)\}$. In the first case, observe that $k\geq \frac{n-1}{2}$ implies $\frac{n-1}{2} \geq n-k-1$. Hence, by relabelling the $x_i$, we may assume that $h[D]> \{h(x_i)\}$ for $1\leq n-k-1$, and so $h$ sends $D'$ to a copy of $S_n^k$, finishing the proof. In the second case, pick an random filter $F''\subseteq P$ which contains all $h(x_i)$ for which $h[D]\bot \{h(x_i)\}$, and which does not contain any element from $h[D]$. Then replacing $h$ by $\turn_{F''}\circ h$ brings us back to the first case.

\textbf{Step 3: From $S_n^k$ ($T_n^k$) to an antichain when $k> \frac{n-1}{2}$}.

We show that if $k> \frac{n-1}{2}$, then any copy of $S_n^k$ in $\P$ can be mapped to an antichain by a permutation in $\G$. Clearly, the dual argument then shows the same for $T_n^k$. Let $\{u_1, \ldots, u_{n-1}\}\subseteq P$ be so that it induces a chain. By the $(n-1)$-transitivity of $\G$, there is some $\rho\in \G$ that maps $\{u_1, \ldots, u_{n-1}\}$ to an antichain $\{v_1, \ldots, v_{n-1}\}$. Let $Z$ be the orbit of $(\P,u_1, \ldots, u_{n-1})$ that is above all the $u_j$. By Lemma~\ref{lem:interface} there exists a canonical function $f: (\P, u_1, \ldots, u_{n-1})\rightarrow (\P, v_1, \ldots, v_{n-1})$ M-generated by $\G$ that agrees with $\rho$ on $\{u_1, \ldots, u_{n-1}\}$. All elements of $Z$ are mapped to one and the same orbit $O$ of $(\P, v_1, \ldots, v_{n-1})$. Now pick $z_1, \ldots, z_k\in Z$ which induce an antichain. By applying an appropriate instance of $\turn$ in a similar fashion as in Step~2, we may assume that $O$ is incomparable with at least $\frac{n-1}{2}$ of the singletons $\{v_i\}$. Choose $(n-k)$ out of these $v_i$. This is possible, as $k> \frac{n-1}{2}$ and consequently $\frac{n-1}{2}\geq n-k$. By relabelling the $u_i$, we may assume that the chosen elements are $v_1, \ldots, v_{n-k}$. Then $f[\{z_1, \ldots, z_k\}]\cup\{v_1, \ldots, v_{n-k}\}$ is an antichain. Since $\{z_1, \ldots, z_k,u_1, \ldots, u_{n-k}\}$ induces a copy of $S_n^k$, we are done.

\textbf{Step 4: From $A_n^k$ to an antichain when $k=\frac{n-1}{2}$}. 

Assuming that $k=\frac{n-1}{2}$, we show that any copy of $A_n^k$ in $\P$ can be mapped to an antichain by a function in $\G$. Note that this assumption implies that $n$ is odd, so $n\geq 5$, and thus $k=\frac{n-1}{2}\geq 2$.

Let $\{x_1, \ldots, x_{k-1}\}\subseteq P$ induce an antichain. Let $s\in P$ be a point below all the $x_i$, and let $\{y_1, \ldots, y_{k}\}\subseteq P$ induce an antichain whose elements are incomparable with all the $x_i$ and $s$. The set $A:=\{s ,x_1, \ldots, x_{k-1}, y_1, \ldots, y_{k}\}$ induces a copy of $A_{n-1}^{k-1}$. By the $(n-1)$-transitivity of $\G$ there exists $\varphi\in \G$ which maps $A$ to an antichain $\{z_1, \ldots, z_{n-1}\}\subseteq P$. Without loss of generality, we write $\varphi(s)=z_{n-1}$, $\varphi(x_i)=z_i$ for $1\leq i\leq k-1$, and $\varphi(y_i)=z_{k+i}$ for $1\leq i\leq k$. By Lemma~\ref{lem:interface} there exists a canonical function $h: (\P, s, x_1, \ldots, x_{k-1}, y_1, \ldots, y_{k})\rightarrow (\P, z_1, \ldots, z_{n-1})$ M-generated by $\G$ which  agrees with $\varphi$ on $A$. Let $U$ be the orbit of $(\P, s, x_1, \ldots, x_{k-1}, y_1, \ldots, y_{k})$ whose elements are larger than $s$ and incomparable to all other elements of $A$. Since $h$ is canonical, $h[U]$ is contained in an orbit $V$ of $(\P, z_1, \ldots, z_{n-1})$. 

Assume that the elements of the orbit $V$ do not satisfy the same relations with all the $z_i$ for $1\leq i\leq n-2$. Then there is a partition $R\cup S=\{z_1, \ldots, z_{n-2}\}$, with both $R$ and $S$ non-empty, such that the elements of $V$ are incomparable with the elements of $R$ and strictly comparable with the elements of $S$. By applying an appropriate instance of $\turn$ we may assume that $|R|\geq k$. Pick any $R'\subseteq R$ of size $k$, any $S'\subseteq S$ of size~1, and a $k$-element antichain $W\subseteq U$. Then  $h^{-1}[R']\cup  h^{-1}[S']\cup W$  induces an antichain of size $n$ whose image $I$ under $h$ induces either $A_n^{k}$ or $B_n^{k}$. In the second case, let $F\subseteq P$ be an random filter which separates the largest element of $I$ from its other elements. Then $\turn_F$ sends $I$ to a copy of $A_n^k$. Thus in either case, $\G$ contains a function which sends an $n$-element antichain to a copy of $A_n^k$. Since $\G$ contains the inverse of all of its functions, it also maps a copy of $A_n^{k}$  to an antichain.

Finally, assume that $V$ satisfies the same relations with all the $z_i$ for $1\leq i\leq n-2$. By applying an appropriate instance of $\turn$ we may assume that $V$ is incomparable with all the $z_i$ for $1\leq i\leq n-2$. Let $W\subseteq U$ induce a $(k-1)$-element antichain, and consider $R:=W\cup \{x_1, y_1, \ldots, y_k, s\}$; then $R$ induces a copy of $A_n^k$. If $V$ is incomparable with $z_{n-1}$, then $h[R]$ is an antichain and we are done. So assume that $V$ and $z_{n-1}$ are comparable. Then $h[R]$ induces $A_n^{k-1}$ or $B_n^{k-1}$. Let $F\subseteq P$ be an random filter that separates $h(s)$ from the other elements of $h[R]$. Then $\turn_F\circ h[R]$ induces $B_n^{n-k+1}$ or $A_n^{n-k+1}$. By Steps~2 and~3, both $A_n^{n-k+1}$ and $B_n^{n-k+1}$ can be mapped to an antichain by permutations from $\G$, finishing the proof.
\end{proof}

\begin{proposition}\label{prop:aboveT}
Let $\G\supsetneq \T$ be a closed group. Then $\G$ contains $\M$.
\end{proposition}
\begin{proof}
If $\G=\Sym_P$, then there is nothing to show, so assume this is not the case. Then $\G$ is not 3-transitive; since $\G\supseteq \T$, the orbits of its action on triples of distinct entries of $P$ are unions of the action of $\T$ on such triples. Since $\G\neq \T$, it cannot preserve $\oz$ or $\od$; thus it preserves $\oo$. Thus $\G\subseteq\M$ by Proposition~\ref{prop:descr:max}. Now if $f\in\G\setminus\T$, then it flips $\oz$ and $\od$, by Lemma~\ref{lem:switch}. Hence, $\rev\circ f$ preserves $\oz$, and so it is an element of $\T\subseteq \G$, by Proposition~\ref{prop:descr:turn}. But then $\rev=\rev\circ f\circ f\inv\in\G$, and so $\G$ contains $\R$.
\end{proof}

Theorem~\ref{thm:groups} now follows from Propositions~\ref{prop:aboveP}, \ref{prop:aboveR}, \ref{prop:aboveM}, \ref{prop:3transitive}, and  ~\ref{prop:aboveT}.

\subsection{Relational description of $\R$}

\begin{prop}\label{prop:descr:rev}
$\Rev=\Aut(P;\bot)$.
\end{prop}
\begin{proof}
	By definition, the function $\rev$ preserves the incomparability relation and its negation, so the inclusion $\supseteq$ is trivial. For the other direction, let $f\in\Aut(P;\bot)$. We claim that $f$ is either an automorphism of $\P$, or satisfies itself the definition of $\rev$ (i.e., $f(b)\leq f(a)$ iff $a\leq b$ for all $a,b\in P$). Suppose that $f$ is not an automorphism of $\P$, and pick $a\leq b$ such that $f(a)\nleq f(b)$. Since $f$ preserves comparability, we then have $f(b)\leq f(a)$. To prove our claim, since $f$ preserves $\bot$ it suffices to show that likewise $f(d)\leq f(c)$ for all $c\leq d$. 
	
	We first observe that if $e\leq b$ and $e\bot a$, then $f(e)\geq f(b)$. For if we had $f(e)\leq f(b)$, then it would follow that $f(e)\leq f(b)\leq f(a)$, a contradiction since $f$ preserves $\bot$. Hence, $f(e)\nleq f(b)$, and so $f(e)\geq f(b)$ since $f$ preserves comparability. 
	
	Next let $r,s\in P$ so that $r\leq s$, $r\leq b$, and $s\bot b$; we show $f(r)\geq f(s)$. Since $f(r)$ and $f(s)$ are comparable, it is enough to rule out $f(r)\leq f(s)$. By our previous observation, we have $f(b)\leq f(r)$, so $f(r)\leq f(s)$ would imply $f(b)\leq f(s)$, contradicting the fact that $f$ preserves $\bot$.
	
	Now  let $u,v\in P$ be so that $u\leq v$ and such that both $u$ and $v$ are incomparable with both $a$ and $b$. Then using the extension property, we can pick $r,s\in P$ as above and such that $u\leq s$ and $v\bot s$. By the preceding paragraph, $f(r)\geq f(s)$, and applying the above once again with $(u,v)$ taking the role of $(r,s)$ and $(r,s)$ the role of $(a,b)$, we conclude $f(v)\geq f(u)$.
		
	Finally, given arbitrary $c,d\in P$ with $c\leq d$, we use the extension property to pick $u,v\in P$ incomparable with all of $a,b,c,d$, and apply the above twice to infer $f(c)\geq f(d)$.
\end{proof}

Theorem~\ref{thm:reducts} now follows from Propositions~\ref{prop:descr:turn}, \ref{prop:descr:max}, and \ref{prop:descr:rev}.

\bibliographystyle{alpha}
\bibliography{randompo_local.bib}

\end{document}